\documentclass{amsart}[11pt]

\usepackage[utf8]{inputenc}
\usepackage{amsmath,amsthm,amsfonts,amssymb,graphicx,amscd}
\usepackage[all,cmtip]{xy}
\usepackage{tikz-cd}
\usetikzlibrary{shapes.geometric}
\usetikzlibrary{calc}
\usepackage{array}
\usepackage{stmaryrd}
\usepackage[english]{babel}
\usepackage{enumerate}
\usepackage{hyperref}
\hypersetup{pdfborder=0 0 0}
\usepackage{tabstackengine}
\usepackage{caption}
\usepackage{subcaption}
\usepackage{float}
\setstackgap{L}{.75\baselineskip}
\usepackage{tensor}
\usepackage{leftindex}
\usepackage{multirow}
\usepackage[normalem]{ulem}

\theoremstyle{plain}

\newsavebox{\myboxa}

\newsavebox{\myboxb}

\newtheorem{theo}{Theorem}
\numberwithin{theo}{section}

\newtheorem{defi}{Def\mbox inition}
\numberwithin{defi}{section}

\numberwithin{conj}{section}

\newtheorem{lemm}{Lemma}
\numberwithin{lemm}{section}

\newtheorem{coro}{Corollary}
\numberwithin{coro}{section}

\numberwithin{prop}{section}

\numberwithin{exam}{section}

\newtheorem{remark}{Remark}
\numberwithin{remark}{section}

\newcommand{\Hom}{\operatorname{Hom}}
\newcommand{\shom}{\operatorname{\underline{Hom}}}
\newcommand{\End}{\operatorname{End}}
\newcommand{\send}{\operatorname{\underline{End}}}
\newcommand{\smod}{\operatorname{\underline{mod}}}
\newcommand{\im}{\operatorname{Im}}
\newcommand{\Ext}{\operatorname{Ext}}
\newcommand{\rad}{\operatorname{rad}}
\newcommand{\soc}{\operatorname{soc}}

\newcommand{\mypicture}[1][]{
\begin{tikzpicture}[#1]
       \tikzset{vnode/.style={draw,thick,circle,minimum width=8mm,inner sep=0pt}};
       \def\dist{2.6}
        
        \begin{scope}
            \node [minimum size=5mm,draw] (zeta0) {{\small $\zeta_0$}};

                \path (zeta0) --++ (504:\dist) node [minimum size=5mm,draw] (zeta2){{\small $\zeta_2$}};
                \path (zeta0) --++ (90:\dist) node [minimum size=5mm,draw] (zeta0) (zeta1) {{\small $\zeta_1$}};
                \path (zeta0) --++ (30:\dist) node [vnode,draw] (zetan+1) {{\small $\zeta_{n+1}$}};
                \path (zeta0) --++ (225:\dist) node  [minimum size=5mm,draw] (zetai-1) {{\small $\zeta_i$}};
                \path (zeta0) --++ (315:\dist) node [vnode,draw] (zetai) {{\small $\zeta_{i+1}$}};
                \path (zeta2) --++ (285:1.3) node (vdots1) {\textbf{$\vdots$}};
                \path (zeta2) --++ (340:3.6) node (vdots2) {\textbf{$\vdots$}};
                \path (zeta2) --++ (140:1.5) node (w) {$\mathcal{W}_{n,\overline{m}}:$};
                
        \end{scope}
        
        \path
            (zeta0) edge node[left] {{\scriptsize $0$}} (zeta1)
                    edge node[left] {{\scriptsize $1$\ \ }} (zeta2)
                    edge node[left] {{\scriptsize $n$\ \ }} (zetan+1)
                    edge node[below] {{\scriptsize $\ \ \ \ \ i-1$}} (zetai-1)
                    edge node[above] {{\scriptsize \ \ \ $i$}} (zetai)
        ;
\end{tikzpicture}
}

\newcommand{\mypictureb}[1][]{
\begin{tikzpicture}[#1]
\foreach \ang\anch in {90/north, 225/south, 135/{north west}}{
  \draw[->,shorten <=7pt, shorten >=7pt] ($(0,0)+(\ang:2)$).. controls +(\ang-40:1.5) and +(\ang+40:1.5) .. ($(0,0)+(\ang:2)$);
}
\def\dist{2.5}
    \node at ($(-5.5,2)+(15:2.5)$) {$Q:$};

    \draw[fill=black] ($(0,0)+(90:2)$) circle (.08);
    \node[anchor=north] at ($(0,-.9)+(90:2.8)$) {$0$};
    \node at ($(0,-1.1)+(90:4.3)$) {$\delta_0$};
    \node at ($(-1.7,-2)+(100:4.2)$) {$\delta_1$};
    \node at ($(-1.5,-6.4)+(100:4.2)$) {$\delta_{i-1}$};

    \draw[fill=black] ($(0,0)+(0:2)$) circle (.08);
    \node[anchor=east] at ($(.5,0)+(0:2.8)$) {$n-1$};

    \draw[fill=black] ($(0,0)+(45:2)$) circle (.08);
    \node[anchor={north east}] at ($(-.1,-.1)+(45:2.8)$) {$n$};
    
    \draw[fill=black] ($(0,0)+(270:2)$) circle (.08);
    \node[anchor=south] at ($(0,.1)+(270:2.8)$) {$i$};

    \draw[fill=black] ($(0,0)+(225:2)$) circle (.08);
    \node[anchor=west] at ($(.8,.5)+(215:2.8)$) {$i-1$};
    
    \draw[fill=black] ($(0,0)+(135:2)$) circle (.08);
    \node[anchor=north west] at ($(.6,-.4)+(135:2.8)$) {$1$};

\draw[->,shorten <=7pt, shorten >=7pt] ($(-.3,4)+(280:2)$) arc (85:140:1.7);
\node at ($(-.4,.2)+(105:2)$) {$\alpha_0$};

\draw[->,shorten <=7pt, shorten >=7pt] ($(0,3)+(228:2.1)$) arc (135:180:2);
\node at ($(-.8,2.5)+(228.5:2.1)$) {$\alpha_1$};

\draw[->,shorten <=7pt, shorten >=7pt] ($(.5,.4)+(300:2)$) arc (325:355:2.8);
\node at ($(.1,1.7)+(338.5:2.5)$) {$\alpha_{n-1}$};
  
\draw[->,shorten <=7pt, shorten >=7pt] ($(-.35,0)+(280:2)$) arc (270:310:2.1);
  \node at ($(.1,4.3)+(315-22.5:2.5)$) {$\alpha_n$};

\draw[->,shorten <=7pt] ($(1,-.6)+(78:2)$) arc (43:80:2.2);
\draw[->,shorten >=7pt] ($(.3,1.3)+(328:2)$) arc (5:50:1.7);
\draw[->,shorten <=7pt] ($(-1.9,2)+(270:2)$) arc (170:210:1.7);
\draw[->,shorten >=7pt] ($(.5,-.5)+(215:2)$) arc (240:270:2.3);
\node at ($(1.7,.2)+(270-20:2.5)$) {$\alpha_i$};
\node at ($(.15,-1.9)+(25:2.5)$) {$\alpha_{n-2}$};
\node at ($(-4.6,-1.5)+(20:2.5)$) {$\alpha_{i-2}\ $};
\node at ($(-1.9,.3)+(315-25:2.5)$) {$\ \alpha_{i-1}$};

\foreach \ang in {310,315,320,176,179.5,183}{
  \draw[fill=black] ($(0,0)+(\ang:2)$) circle (.02);
}
\end{tikzpicture}
}

\begin{document}

\title[UDR of modules over generalized Brauer tree algebras]{Universal deformation rings of a special class of modules over generalized Brauer tree algebras}
\date{}
\author[Caranguay-Mainguez]{Jhony F. Caranguay-Mainguez}
\address{Universidad de Antioquia, Instituto de Matem\'aticas}
\email{jhony.caranguay@udea.edu.co}

\author[Rizzo]{Pedro Rizzo}
\address{Universidad de Antioquia, Instituto de Matem\'aticas}
\email{pedro.hernandez@udea.edu.co}

\author[V\'elez-Marulanda]{Jos\'e A. V\'elez-Marulanda}
\address{Department of Applied Mathematics \& Physics, Valdosta State University, Valdosta, GA,  USA}
\email{javelezmarulanda@valdosta.edu}
\address{Facultad de Matem\'aticas e Ingenier\'{\i}as, Fundaci\'on Universitaria Konrad Lorenz, Bogot\'a D.C.,  Colombia}
\email{josea.velezm@konradlorenz.edu.co}

\maketitle

\begin{abstract}
Let $\Bbbk$ be an algebraically closed f\mbox ield and $\Lambda$ a generalized Brauer tree algebra over $\Bbbk$. We compute the universal deformation rings of the periodic string modules over $\Lambda$. Moreover, for a specif\mbox ic class of generalized Brauer tree algebras $\Lambda(n,\overline{m})$, we classify the universal deformation rings of the modules lying in $\Omega$-stable components $\mathfrak{C}$ of the stable Auslander-Reiten quiver provided that $\mathfrak{C}$ contains at least one simple module. Our approach uses several tools and techniques from the representation theory of Brauer graph algebras. Notably, we leverage Duf\mbox f\mbox ield's work on the Auslander-Reiten theory of these algebras and Opper-Zvonareva's results on derived equivalences between Brauer graph algebras.
\end{abstract}

\section{Introduction}

Let $\Bbbk$ be an algebraically closed f\mbox ield, $\widehat{\mathcal{C}}$ be the category of complete local Noetherian commutative $\Bbbk$\hspace{-1mm} -\hspace{-1mm} algebras with residue f\mbox ield $\Bbbk$, and $\Lambda$ be a $\Bbbk$\hspace{-1mm} -\hspace{-1mm} algebra of f\mbox inite dimension.

Deformation theory studies lifts of a f\mbox initely generated $\Lambda$\hspace{-1mm} -\hspace{-1mm} module $M$ to a f\mbox initely generated $R\otimes_{\Bbbk} \Lambda$\hspace{-1mm} -\hspace{-1mm} module $W$, where $R\in Ob(\widehat{\mathcal{C}})$ and $W$ is free over $R$. This theory, inspired by Mazur's deformation theory for f\mbox inite-dimensional representations of prof\mbox inite groups, has become a signif\mbox icant tool in number theory and arithmetic geometry (see \cite{bc5}). Universal deformation rings, which parametrize isomorphism classes of lifts of a Galois representation, play a crucial role in this context.

The scope of deformation theory has expanded beyond its original domain. It has been successfully extended to modules over f\mbox inite-dimensional algebras (e.g., \cite{bv1,ile,lau,yau}) and, more recently, inf\mbox inite-dimensional algebras (e.g., \cite{adr,lop}). In particular, F. M. Bleher and the third author showed in \cite{bv1} that every module over a self-injective algebra with stable endomorphism ring isomorphic to $\Bbbk$ admits a universal deformation ring. Furthermore, for Frobenius algebras, they showed that universal deformation rings are preserved under the syzygy functor $\Omega$. Therefore, a central problem is the classif\mbox ication of universal deformation rings for $\Lambda$-modules $M$ satisfying the condition $\send_{\Lambda}(M)\cong \Bbbk$. Developing classif\mbox ication results for these rings, along with deformation theory for $\Lambda$, provides tools and novel invariants to clarify the intricate structure of the category mod-$\Lambda$ and improve our understanding of its objects.

The results mentioned above enable the application of representation theory techniques to compute universal deformation rings. Special biserial algebras, a signif\mbox icant class of algebras with connections to group algebras, algebraic geometry, and derived categories (see \cite{syb}), play a crucial role in representation theory. The combinatorial nature of special biserial algebras allows for a detailed understanding of their representation theory. Butler and Ringel's work in \cite{br}, building on Gelfand and Ponomarev's results in \cite{gp}, shows that non-projective indecomposable modules are string and band modules, and irreducible homomorphisms arise from adding or removing hooks and cohooks. Moreover, Erdmann and Skowro\'nski classif\mbox ied the shapes of components in the stable Auslander-Reiten quiver of self-injective special biserial algebras in \cite{kes}. Brauer graph algebras, a subclass of special biserial algebras, have been studied extensively. For instance, Duf\mbox f\mbox ield's work in \cite{duff} provides criteria for the location of modules in components of the stable Auslander-Reiten quiver. These algebras originate from the modular representation theory of f\mbox inite groups and appear in various contexts within representation theory. Notably, all symmetric special biserial algebras are Brauer graph algebras, and every special biserial algebra is a quotient of a symmetric special biserial algebra (see \cite{sib, syb}).

Since symmetric algebras are Frobenius, we can apply Bleher-Vélez's deformation theory for modules over Frobenius algebras, combined with the representation theory of special biserial algebras developed in \cite{br}, \cite{kes}, and \cite{duff}, to study deformation of f\mbox inite dimensional modules over Brauer graph algebras. Research on universal deformation rings of modules over Brauer graph algebras has been directed, particularly for \textit{generalized Brauer trees} (Brauer graphs without cycles). In \cite{bw}, Bleher and Wackwitz determined the universal deformation rings of f\mbox initely generated modules over certain Nakayama algebras, which includes f\mbox inite-representation type Brauer tree algebras. Meyer et al. extended this work to inf\mbox inite-representation type generalized Brauer tree algebras of polynomial growth in \cite{msw}. A key tool in these studies is derived equivalence. In particular, Kauer moves and Opper-Zvonareva's classif\mbox ication of derived equivalences (see \cite{oz}) allow for a reduction to simpler cases.

Other important examples of research universal deformation rings of modules over Brauer graph algebras are given in \cite[Section 3]{bv1}, \cite{vm} and \cite{cgrv}. While the algebras studied in \cite{bv1} and \cite{vm} are implicitly generalized Brauer tree algebras (being associated with stars), we extend these results to a broader class of generalized Brauer tree algebras for the case of periodic string modules and certain non-periodic string modules.

We compute universal deformation rings of periodic string modules $M$ over generalized Brauer tree algebras with stable endomorphism ring isomorphic to $\Bbbk$. Our approach utilizes derived equivalences involving stars, leveraging Opper-Zvonareva's classif\mbox ication of derived equivalences for Brauer graph algebras. As an application, we determine the universal deformation rings of periodic string modules over standard Koszul symmetric special biserial algebras, relying on Magyar's classif\mbox ication in \cite[Prop. 3.8]{mgy}. Finally, we examine the $\Omega$-stable components of type $\mathbb{Z}A_{\infty}^{\infty}$ containing simple modules over a generalized Brauer tree algebra given by a star.

This paper is organized as follows. In Section \ref{cprel} we provide background on deformation theory of modules over self-injective algebras and the representation theory of Brauer graph algebras, including derived equivalences and the structure of their Auslander-Reiten quivers. In Section \ref{sec:main}, we identify periodic string modules over generalized Brauer tree algebras with one-dimensional stable endomorphism rings, compute their tangent spaces and universal deformation rings, and apply these results to standard Koszul algebras. Finally, in Section \ref{nonperiodiccase}, for a class of generalized Brauer tree algebras $\Lambda(n,\overline{m})$ of non-polynomial growth (see Subsection \ref{sedis} for more details about the def\mbox inition of these algebras), we extend our analysis to the modules in $\Omega$-stable components $\mathfrak{C}$ of the stable Auslander-Reiten quiver of $\Lambda(n,\overline{m})$ when $\mathfrak{C}$ contains at least one non-periodic simple module. We compute and classify their universal deformation rings.

This article constitutes the main results of the f\mbox irst author's doctoral dissertation, supervised by the other two authors.

\section{Preliminaries}\label{cprel}

This section provides the essential background to facilitate understanding of our main results and the computations involved in the following sections. This foundational material includes a brief introduction to deformation theory of modules over algebras, a concise overview of Brauer graph algebras, a study of their derived equivalences, and certain features of their Auslander-Reiten quivers. Within Brauer graph algebras, we will focus on one particular type, known as \textit{generalized Brauer tree algebras} of non-polynomial growth.

We begin by f\mbox ixing some notations used throughout the paper. The ring $\Bbbk\llbracket x_1,\ldots,x_n\rrbracket$ represents the power series ring in $n$ variables with coef\mbox f\mbox icients in $\Bbbk$, with $\Bbbk$ denoting an algebraically closed f\mbox ield of arbitrary characteristic. By an algebra, we mean a f\mbox inite-dimensional unital associative $\Bbbk$-algebra. For every algebra $\Lambda$, by a $\Lambda$-module $M$ we mean a f\mbox initely generated right $\Lambda$-module unless otherwise specif\mbox ied. 

For any pair $M$, $N$ of $\Lambda$-modules we denote the \textit{space of stable homomorphisms} by $\shom_{\Lambda}(M,N)$, which is def\mbox ined to be the quotient $\Hom_{\Lambda}(M,N)/\mathcal{P}_{\Lambda}(M,N)$, where $\mathcal{P}_{\Lambda}(M,N)$ is the $\Bbbk$-vector subspace of $\Hom_{\Lambda}(M,N)$ generated by all the module homomorphisms that factor through projective modules. In particular, we denote by $\End_{\Lambda}(M)=\Hom_{\Lambda}(M,M)$ the endomorphism ring of $M$ and by $\send_{\Lambda}(M)$ the \textit{stable endomorphism ring} of $M$. We denote by\hspace{-2mm} $\mod\hspace{-1.5mm} -\hspace{-.3mm} \Lambda$ the category of $\Lambda$-modules and write $\smod\hspace{-.4mm}-\Lambda$ to denote the stable category of\hspace{-2mm} $\mod\hspace{-1.5mm} -\hspace{-.3mm} \Lambda$. Recall that the objects of the category $\smod-\Lambda$ are the same as for\hspace{-2mm} $\mod\hspace{-1.5mm} -\hspace{-.3mm} \Lambda$ and as morphisms we have $\shom_{\Lambda}(M,N)$, for any $M$ and $N$ $\Lambda$-modules. 

In the category of $\Lambda$-modules, we introduce two key operators: the \textit{syzygy of} $M$, denoted by $\Omega(M)$, and the \textit{Auslander-Reiten translation}, denoted by $\tau$. The syzygy is (an endofunctor) def\mbox ined as the kernel of a projective cover of $M$, while the Auslander-Reiten translation is an endofunctor that maps modules their ``neighbors'' in the Auslander-Reiten quiver. For more details, we recommend \cite{sy}.

Finally, for any two morphisms $f:A\rightarrow B$ and $g:B\rightarrow C$ in a category $\mathcal{C}$, we write $gf$ to denote the composition. Moreover, for every f\mbox inite set $S$, we write $|S|$ to denote the number of elements of $S$.

\subsection{Deformation theory: An overview}

We will denote by $\widehat{\mathcal{C}}$ the category whose objects are commutative complete local Noetherian $\Bbbk$-algebras with residue f\mbox ield $\Bbbk$. Let $\Lambda$ be an algebra, $M$ be a $\Lambda$-module and $R\in Ob(\widehat{\mathcal{C}})$. 

A \emph{lift of $M$ over $R$} is a pair $(W,\phi)$, where $W$ is a f\mbox initely generated right $\Lambda\otimes_{\Bbbk}R$-module which is free over $R$, and $\phi: W\otimes_R \Bbbk\rightarrow M$ is an isomorphism of $\Lambda$-modules. We say that two lifts $(W,\phi)$ and $(W',\phi')$ of $M$ over $R$ are \textit{isomorphic} if there is an isomorphism $f:W\rightarrow W'$ of $R\otimes_{\Bbbk} \Lambda$-modules such that $\phi=\phi'(id_{\Bbbk}\otimes f)$, where $id_{\Bbbk}$ is the identity map over $\Bbbk$. A \emph{deformation of $M$ over} $R$ is an isomorphism class of a lift of $M$ over $R$. We denote by $\text{Def}_{\Lambda}(M,R)$ the set of all deformations of $M$ over $R$. This set has a functorial behavior, in the following sense. We def\mbox ine the covariant \emph{deformation functor} $\widehat{F}_M:\widehat{\mathcal{C}}\rightarrow\textbf{Sets}$ as follows: for any $R\in \text{Obj}(\widehat{\mathcal{C}})$, we put $\widehat{F}_M(R):=\text{Def}_{\Lambda}(M,R)$ and for any morphism $\alpha:R\rightarrow R'$ in $\widehat{\mathcal{C}}$, the function $\widehat{F}_M(\alpha):\text{Def}_{\Lambda}(M,R)\rightarrow \text{Def}_{\Lambda}(M,R')$ is def\mbox ined by $\widehat{F}_M(\alpha)([(W,\phi)])=[(W\otimes_{R,\alpha}R',\phi_{\alpha})]$, for every $[(W,\phi))]\in \text{Def}_{\Lambda}(M,R)$, where $\phi_{\alpha}:(W\otimes_{R,\alpha}R')\otimes_{R'}\Bbbk\rightarrow M$ is obtained by the composition of isomorphisms $(W\otimes_{R,\alpha}R')\otimes_{R'}\Bbbk\simeq W\otimes_{R}\Bbbk\stackrel{\phi}{\rightarrow} M$. If there exists a unique object $R(\Lambda,M)$ in $\widehat{\mathcal{C}}$ such that $\widehat{F}_M$ is naturally isomorphic to the functor $\text{Hom}_{\widehat{\mathcal{C}}}(R(\Lambda,M),-)$, i.e., $R(\Lambda,M)$ represents $\widehat{F}_M$, we assert that $R(\Lambda,M)$ is the \emph{universal deformation ring} of $M$. An important consequence arising from this situation is that there exists a deformation $[(U(\Lambda, M),\phi_{U(\Lambda, M)})]$, which we call the \emph{universal deformation of $M$}, such that, given any $R\in \text{Obj}(\widehat{\mathcal{C}})$ and any $[(W,\phi)]\in Def_{\Lambda}(M,R)$, there exists a unique morphism $\alpha:R(\Lambda, M)\rightarrow R$ such that $\widehat{F}_M(\alpha)([\left(U(\Lambda, M),\phi_{U(\Lambda, M)}\right)])=[(W,\phi)]$.

Moreover, the \textit{tangent space} $t_M$ \textit{of} $F_M$ is def\mbox ined by $t_M=F_M(\Bbbk\llbracket x \rrbracket/\langle x^2\rangle)$. As consequence from \cite[Lemma 2.10]{schl} we obtain that $t_M$ is a vector space, which is isomorphic to the f\mbox irst group of extensions $\Ext_{\Lambda}^1(M,M)$ (see \cite[Proposition 2.1]{bv1}). Besides, if $M$ admits a universal deformation ring then, $R(\Lambda,M)=\Bbbk$ if $t_M=0$ and $R(\Lambda,M)$ is a quotient of $\Bbbk \llbracket x_1,\ldots,x_{\dim_{\Bbbk}(t_M)}\rrbracket$ if $\dim_{\Bbbk}(t_M)>0$ (see \cite[pg. 223]{bw}).

The following theorem, which is an adapted version to our context of Theorem 1.1 in \cite{rv}, will be useful to compute the universal deformation ring in several subsequent situations.

\begin{theo}[\cite{rv}]\label{anote} Let $\Lambda$ be a $\Bbbk$-algebra and $M$ be an indecomposable right $\Lambda$-module with $\dim_{\Bbbk} M<\infty$ and $\dim_{\Bbbk} \Ext_{\Lambda}^1(M,M)=1$. Assume that $M$ has a universal deformation ring $R(\Lambda,M)$, and that there exists an ordered sequence of indecomposable f\mbox inite dimensional right $\Lambda$-modules (up to isomorphism) $\mathcal{L}_M=\{W_0,W_1,\ldots\}$, with $W_0=M$ and such that for $l\geq1$, there are a monomorphism $\iota_l:W_{l-1}\rightarrow W_l$ and an epimorphism $\epsilon:W_l\rightarrow W_{l-1}$ such that the composition $\sigma_l:=\iota_l\epsilon_l$ satisf\mbox ies that $\ker(\sigma_l)\cong M$, $\im(\sigma_l^l)\cong M$ and $\mathcal{L}_M$ is maximal in the sense that if there is another ordered sequence of indecomposable right $\Lambda$-modules $\mathcal{L}_{M'}$ with these properties, then $\mathcal{L}_{M'}\subseteq \mathcal{L}_M$.
\begin{enumerate}
\item[$i)$] If $\mathcal{L}_M$ is f\mbox inite, and its last element, say $W_N$, satisf\mbox ies the conditions $\dim_{\Bbbk} \Hom_{\Lambda}(W_N,M)=1$ and $\Ext_{\Lambda}^1(W_N,M)=0$, then $R(\Lambda,M)\cong \Bbbk \llbracket x \rrbracket/\langle x^{N+1}\rangle$.
\item[$ii)$] If $\mathcal{L}_M$ is inf\mbox inite, then $R(\Lambda,M)\cong \Bbbk \llbracket x \rrbracket$.
\end{enumerate}
\end{theo}

Recall that $\Lambda$ is \emph{self-injective algebra} if the right $\Lambda$-module $\Lambda_{\Lambda}$ is injective. On the other hand, we say that $\Lambda$ is a \emph{Frobenius algebra} if there exists an isomorphism as right $\Lambda$-modules between $\Lambda_{\Lambda}$ and $D(\Lambda)_{\Lambda}$, where $D(\Lambda)=\text{Hom}_{\Bbbk}(\Lambda_{\Lambda},\Bbbk)$. Recall also that $\Lambda$ is called \emph{symmetric algebra} if $\Lambda$ is endowed with a non-degenerate, associative bilinear form $B:\Lambda\times\Lambda\rightarrow\Bbbk$ such that $B(a,b)=B(b,a)$, for all $a,b\in\Lambda$. By \cite[Chapter IV, Proposition 3.8]{sy}, any Frobenius algebra is self-injective and every symmetric algebra is Frobenius.

The central aim of deformation theory is to determine when an indecomposable module admits a universal deformation ring and to calculate it. In order to resolve these problems the following theorem will be a useful tool.

\begin{theo}[\cite{bv1}]\label{trsyz} Suppose that $\Lambda$ is a self-injective algebra and let $M$ be an indecomposable $\Lambda$-module. The following statements hold.
\begin{enumerate}
   \item[$i)$] If $\send_{\Lambda}(M)\cong \Bbbk$, then $M$ admits a universal deformation ring, $R(\Lambda,M)$.
    \item[$ii)$] If $\Lambda$ is a Frobenius algebra and $\send_{\Lambda}(M)\cong \Bbbk$, then $\send_{\Lambda}(\Omega(M))\cong \Bbbk$ and $R(\Lambda,\Omega(M))\cong R(\Lambda,M)$ in $\widehat{\mathcal{C}}$.
\end{enumerate}
\end{theo}

Also, we are interested in computing the dimension of $\Ext_{\Lambda}^1(M,M)$ for a module $M$ over an algebra $\Lambda$. If $\Lambda$ is a self-injective algebra and $M$ is a $\Lambda$\hspace{-1mm} -\hspace{-1mm} module, then $\Ext_{\Lambda}^1(M,M)\cong \shom_{\Lambda}(\Omega(M),M)$ (see \cite[Lemma 5.1]{sk}).

\subsection{Brauer Graph algebras: A short background}

Brauer graph algebras are well-known to coincide with symmetric special biserial algebras (see, e.g., \cite{sib}), a class of algebras with tame representation type (see, e.g., \cite{syb}). Their signif\mbox icance in representation theory stems from the fact that every special biserial algebra can be expressed as a quotient of a symmetric special biserial algebra. Recall that our focus will lie on generalized Brauer tree algebras of non-polynomial growth.

Next, for the sake of completeness, we will introduce the def\mbox initions of special biserial and Brauer graph algebras.

\begin{defi}[\cite{syb}]\label{dsb} Let $Q$ be a quiver and $I$ an admisible ideal of $\Bbbk Q$. We say that the pair $(Q,I)$ is special biserial if the following conditions hold.
\begin{enumerate}
\item[(G1)] Every vertex is the starting point of at most two arrows and is the ending point of at most two arrows.

\item[(G2)] For every arrow $\beta$, there is at most one arrow $\alpha$ such that $t(\alpha)=s(\beta)$ and $\alpha \beta \notin I$, and there is at most one arrow $\gamma$ such that $t(\beta)=s(\gamma)$ and $\beta \gamma\notin I$.
\end{enumerate}
\end{defi}

Now, we def\mbox ine special biserial algebras and a particular class of special biserial algebras, the so-called \textit{string algebras}, which will be important in our calculations.

\begin{defi}[\cite{syb,br}]\label{dsb} Let $\Lambda$ be an algebra.
\begin{enumerate}[$i)$]
    \item $\Lambda$ is called \emph{special biserial algebra} if it is Morita-equivalent to a bound quiver algebra $\Bbbk Q/I$ such that $(Q,I)$ is a special biserial pair.
    \item $\Lambda$ is called \emph{string algebra} if it is Morita-equivalent to a bound quiver algebra $\Bbbk Q/I$ such that $(Q,I)$ is a special biserial pair and $I$ is generated by zero relations.
\end{enumerate}
\end{defi}

To give the precise def\mbox inition of Brauer graph algebras, we f\mbox irst present the def\mbox inition of \textit{ribbon graph} following \cite[Section 2.7]{syb} and \cite[Def\mbox inition 1.1]{oz}.

\begin{defi}[\cite{syb},\cite{oz}]\label{dribbg} A \emph{ribbon graph} is a tuple $\Gamma=(V,H,s,\iota,\sigma)$ consisting of the following data.
\begin{enumerate}[$i)$]
    \item $V$ is a f\mbox inite set whose elements are called \emph{vertices};
    \item $H$ is a f\mbox inite set whose elements are called \emph{half-edges};
    \item $s:H\rightarrow V$ is a function;
    \item $\iota:H\rightarrow H$ is a function without f\mbox ixed points;
    \item $\sigma:H\rightarrow H$ is a permutation such that the orbits in $\sigma$ correspond to the sets $s^{-1}(v)$ for $v\in V$.
\end{enumerate}
\end{defi}

Following \cite[p.7]{oz}, for a ribbon graph $\Gamma=(V,H,s,\iota,\sigma)$, we def\mbox ine an \textit{edge} as a set of the form $\{h,\iota(h)\}$. From now on, we use symbols like $i$ and $j$ to denote edges and we use $i_+$ and $i_-$ to denote the half-edges corresponding to an edge $i$, and we def\mbox ine $\overline{i_-}=i_+$ and $\overline{i_+}=i_-$. So, in this case, $i=\{i_+,i_-\}$.

\begin{defi}[\cite{syb,oz}]\label{dbra} A \emph{Brauer graph} is a quadruple $\mathcal{G}=(V,E,\mathfrak{m},\mathfrak{o})$ consisting of the following data:

\begin{enumerate}[$i)$]
    \item A pair $(V,E)$, which is a f\mbox inite and connected graph (loops and multiple edges are allowed) such that $E$ is the set of edges corresponding to a ribbon graph of the form $\Gamma=(V,H,s,\iota,\sigma)$.
    \item For each $v\in V$, we def\mbox ine the \emph{valency} of $v$, denoted by $val(v)$, as the number of half-edges incident to $v$.
    \item A function $m:V\rightarrow \mathbb{Z}^+$ called multiplicity function. We say that a vertex $v\in V$ is \emph{truncated} if $val(v)m(v)=1$. We denote by $V^{\ast}$ the set of non-truncated vertices of $\mathcal{G}$.
  \item A function $\mathfrak{o}$, called \emph{orientation}, that assigns to each non-truncated vertex $v\in V$ a cyclic order induced by the $\sigma$-orbit on the edges incident to $v$ such that, if $val(v)=1$, then $\mathfrak{o}(v)$ is given by $i<i$, where $i$ is the unique incident edge to $v$. Thus, for every non-truncated vertex $v\in V$, the orientation $\mathfrak{o}(v)$ can be written as $i_0<i_1<\cdots<i_{val(v)-1}<i_0$, where $i_0,i_1,\cdots,i_{val(v)-1}$ are all the incident edges to $v$. In this case, we say that $i_{k+1}$ is a \emph{successor} of $i_k$ if $0\leq k<val(v)-1$ and $i_0$ is the successor of $i_{val(v)-1}$.
\end{enumerate}
\end{defi}

We represent any Brauer graph $\mathcal{G}$ into an oriented plane such that $\mathfrak{o}$ is given by a counterclockwise orientation (we do not assume that the edges do not cross each other). For every $v\in V$ and every edge $i$ incident to $v$, the cyclic order $\mathfrak{o}$ has the form $i_0:=i<i_1<\cdots<i_{val(v)-1}<i_0$. We call to $i_0,\ldots,i_{val(v)-1}$ a \textit{successor sequence} of $v$ starting at $i$. An edge can be appear twice in the cyclic order of a vertex $v$ and hence, in general, there is not a unique successor sequence of $v$ starting at $i$. For this reason, we must consider half-edges. If $i\in E$ is not a loop, then we just write $i_+=i_-=i$. Observe that every pair consisting of a vertex $v$ and a half-edge $i_{\ast}$ of $i\in E$ incident to $v$, with $\ast\in \{+,-\}$, determines a unique successor sequence of $v$ starting at $i$. We denote by $\zeta(v,i_{\ast})$ to such successor sequence.

\begin{defi}[\cite{syb}]\label{dbga} Let $\mathcal{G}=(V,E,\mathfrak{m},\mathfrak{o})$ be a Brauer graph. We def\mbox ine the \emph{Brauer graph algebra} associated to $\mathcal{G}$ as the bound quiver algebra $\Lambda_{\mathcal{G}}:=\Bbbk Q_{\mathcal{G}}/I_{\mathcal{G}}$, where the quiver $Q_{\mathcal{G}}=(Q_0,Q_1,s,t)$ is given as follows:
\begin{itemize}
    \item[$i)$] $Q_0=E$.
    \item[$ii)$] Let $i,j\in Q_0$. We put an arrow $\alpha\in Q_1$ from $i$ to $j$ provided that $i$ and $j$ are incident to a common non-truncated vertex $v\in V$ and $j$ is a successor of $i$ in the cyclic order of $v$. Note that if $v$ is a non-truncated vertex and $i$ is an edge incident to $v$, then, for each half-edge $i_{\ast}$ of $i$, with $\ast\in \{+,-\}$, the successor sequence $\zeta(v,i_{\ast})$ induces a cyclic path (starting at $i$) $\alpha_0\cdots \alpha_{val(v)-1}$ in $\mathcal{P}(Q_{\mathcal{G}})$, with $\alpha_k\in Q_1$ for $0\leq k\leq val(v)-1$. This path is called a \emph{special} $v$-\emph{cycle} and is denoted by $A_{v,i_{\ast}}$.
    \item[$iii)$] The ideal $I_{\mathcal{G}}$ of the path algebra $\Bbbk Q_{\mathcal{G}}$ is generated by the relations of type I, II and III def\mbox ined as follows.
    \begin{enumerate}
        \item[{\bf Type I:}] $A_{v,i_+}^{\mathfrak{m}(v)}-A_{v',i_-}^{\mathfrak{m}(v')}$, where $v$ and $v'$ are non-truncated vertices and $i$ is an edge incident to $v$ and $v'$.
        \item[{\bf Type II:}] $A_{v,i_{\ast}}^{\mathfrak{m}(v)}A_{v,i_{\ast}}^0$, where $v$ is a non-truncated vertex, $i$ is an edge incident to $v$ and $\ast \in \{+,-\}$. Here, $A_{v,i_{\ast}}^0$ denotes the f\mbox irst arrow of $A_{v,i_{\ast}}$.
        \item[{\bf Type III:}] $\alpha \beta$, where $\alpha, \beta \in Q_1$, $\alpha \beta$ is not a subpath of any special cycle, except if $\alpha=\beta$ is a loop corresponding to a non-truncated vertex $v\in V$ such that $val(v)=1$.
    \end{enumerate}
\end{itemize}
\end{defi}

The following theorem was proved by S. Schroll in \cite{sib} and constitutes an important result on representation theory of special biserial algebras.

\begin{theo}[\cite{sib}]\label{tbgss} The class of Brauer graph algebras coincides with the class of symmetric special biserial algebras.
\end{theo}

\begin{defi}[\cite{syb}]\label{dbta} A \emph{generalized Brauer tree algebra} is a Brauer graph algebra $\Lambda_{\mathcal{G}}$, with $\mathcal{G}=(V,E,m,\mathfrak{o})$, such that the graph $(V,E)$ is a tree, i.e., it does not contain cycles. Moreover, we say that a generalized Brauer tree algebra $\Lambda_{\mathcal{G}}$ is called a Brauer tree algebra if there exists at most one vertex in $V$ with multiplicity greater than one.
\end{defi}

\subsection{Derived equivalences between Brauer graph algebras}

A key tool in our calculations is the derived equivalence for Brauer graph algebras. Our interest in derived equivalences is motivated by the results in \cite{bv2} where it is shown that universal deformation rings are preserved under stable equivalences of Morita type. We now introduce the essential background for these results.

Two algebras $\Lambda_1$ and $\Lambda_2$ are considered \textit{derived equivalent} if their derived categories $D^b(\hspace{-2.6mm} \mod\hspace{-1.6mm} -\hspace{-.5mm} \Lambda_1)$ and $D^b(\hspace{-2.6mm} \mod\hspace{-1.6mm} -\hspace{-.5mm} \Lambda_2)$ are equivalents as triangulated categories (see \cite[Section 9.2]{kr21}). Within the framework of Brauer graph algebras, the seminal work on derived equivalences is attributed to Antipov et al. (see \cite{az}, \cite{oz}). To establish the necessary groundwork for their f\mbox indings, we introduce the following preliminaries.

Let $\mathcal{G}=(V,E,\mathfrak{m},\mathfrak{o})$ a Brauer graph algebra with underlying ribbon graph $(V,H,s,\iota,\sigma)$. Let $i$ and $j$ be edges and $\ast,\star \in \{+,-\}$. Assume that $s(i_{\ast})=s(j_{\star})=v$ for a vertex $v\in V$. If $j$ is the successor of $i$ in a successor sequence induced by the positions of $i_{\ast}$ and $j_{\star}$, then we write $j_{\ast}=\sigma(i_{\star})$.

\begin{defi}[\cite{oz}]\label{dface} Let $\mathcal{G}$ be a Brauer graph with associated Brauer graph algebra $\Lambda$ and $p$ be a positive integer. A face of $\mathcal{G}$ of perimeter $p$ is an equivalence class of a cyclic sequence of half-edges
$$
F=((h_1)_{\ast_1},\ldots,(h_{2p})_{\ast_{2p}}),
$$
where $\ast_t\in \{+,-\}$ for all $t\in \{1,\cdots, 2p\}$, $(h_{2k+2})_{\ast_{2k+2}}=\overline{(h_{2k+1})_{\ast_{2k+1}}}$ and $(h_{2k+1})_{\ast_{2k+1}}=\sigma((h_{2k})_{\ast_{2k}})$, $F$ is not a power of another sequence with these conditions, and the equivalence relation is given by cyclic permutations.
\end{defi}

\begin{theo}[\cite{oz}]\label{tzo} Let $\Lambda_1$ and $\Lambda_2$ be Brauer graph algebras with associated Brauer graphs $\mathcal{G}_1$ and $\mathcal{G}_2$. Then, $\Lambda_1$ and $\Lambda_2$ are derived equivalent if and only if all the following conditions hold.
\begin{enumerate}
    \item[$i)$] $\mathcal{G}_1$ and $\mathcal{G}_2$ have the same number of vertices, the same number of edges, and the same number of faces;
    \item[$ii)$] $\mathcal{G}_1$ and $\mathcal{G}_2$ have the same multisets of multiplicities and the same multiset of perimeters;
    \item[$iii)$] $\mathcal{G}_1$ and $\mathcal{G}_2$ are both bipartite or neither of them is a bipartite graph.
\end{enumerate}
\end{theo}

Note that Theorem \ref{tzo} employs the multisets of multiplicities and perimeters associated with a Brauer graph. Specif\mbox ically, for a Brauer graph $\mathcal{G}=(V,E,\mathfrak{m},\mathfrak{o})$ the multiset of multiplicities is def\mbox ined as $\{\mathfrak{m}(v):v\in V\}$, and the multiset of perimeters is given by $\{p(F):F \in F(\mathcal{G})\}$, where $F(\mathcal{G})$ represents the set of faces of $\mathcal{G}$ and $p(F)$ denotes the perimeter of a face $F$ in $F(\mathcal{G})$.

An important consequence of derived equivalences in this work is their impact on equivalences between deformation theories. Notably, derived equivalences between self-injective algebras induce stable equivalences of Morita type (Rickard, \cite[Corollary 5.5]{jr}). Furthermore, universal deformation rings are invariant under such equivalences (\cite[Proposition 3.2.6]{bv2}). We now present these results in detail.

\begin{defi}[\cite{bv2,brr}]\label{dsemt} We say that two f\mbox inite-dimensional algebras $\Lambda_1$ and $\Lambda_2$ are \emph{stably equivalent of Morita type} if there exist a $\Lambda_2-\Lambda_1$-bimodule $X$ and a $\Lambda_1-\Lambda_2$-bimodule $Y$, both projective as left and right modules, and
$$
X\otimes_{\Lambda_1}Y\cong \Lambda_2 \oplus Q,
$$
$$
Y\otimes_{\Lambda_2}X\cong \Lambda_1 \oplus P,
$$
where $P$ is a projective $\Lambda_1-\Lambda_1$-bimodule and $Q$ is a projective $\Lambda_2-\Lambda_2$-bimodule.
\end{defi}

Stable equivalences of Morita type, as introduced by M. Broué in \cite[p.19]{brr} (see also \cite[Def\mbox inition 3.2.3]{bv2}), are crucial for investigating invariants of stable categories. Specif\mbox ically, a stable equivalence of Morita type between algebras $\Lambda_1$ and $\Lambda_2$, induces an equivalence between their corresponding stable module categories $\smod-\Lambda_1$ and $\smod-\Lambda_2$.

\begin{theo}[\cite{jr}]\label{tderst} Let $\Lambda_1$ and $\Lambda_2$ be two self-injective and f\mbox inite-dimensional $\Bbbk$-algebras which are derived equivalent. Then, $\Lambda_1$ and $\Lambda_2$ are stably equivalent of Morita type.
\end{theo}

\begin{theo}[\cite{bv2}]\label{tedstm} With the hypothesis in Def\mbox inition \ref{dsemt}, if $M$ is a $\Lambda_1$-module and $M'$ is the $\Lambda_2$-module $X\otimes_{\Lambda_1} M$, then the deformation functors $\widehat{F}_M$ and $\widehat{F}_{M'}$ are naturally isomorphic. In particular, $M$ admits a universal deformation ring if and only if $M'$ admits a universal deformation ring. In this situation, the universal deformation rings $R(\Lambda_1,M)$ and $R(\Lambda_2,M')$ are isomorphic in $\mathcal{\widehat{C}}$.
\end{theo}

\subsection{A brief study of Auslander-Reiten quiver properties}

\subsubsection{The Auslander-Reiten quiver of a Brauer graph algebra}

Recall that, by \cite[I.8.11. p.26]{kebt}, the \textit{stable Auslander-Reiten quiver} $\leftindex_s\Gamma_{\Lambda}$ of an algebra $\Lambda$ coincides with the Auslander-Reiten quiver of the algebra $\Lambda_s:=\Lambda/S$, where $S$ is the direct sum of all the socles of projective-injective indecomposable modules over $\Lambda$. In particular, if $\Lambda$ is a self-injective algebra, then $\Lambda_s=\Lambda/\soc(\Lambda)$.

Now, if $\Lambda$ is a special biserial algebra, then $\Lambda_s$ is a string algebra, by \cite[II.1.3. p.47]{kebt}. Hence, we are interested in study the components of the Auslander-Reiten quiver of a string algebra. Moreover, recall that if $\mathfrak{C}$ is a component of $\leftindex_s\Gamma_{\Lambda}$, we say that $\mathfrak{C}$ is an \textit{exceptional component} provided is composed entirely by string modules; we say that $\mathfrak{C}$ is an \textit{homogeneous component} if it is composed entirely by band modules.

The following result given by K. Erdmann and A. Skowro\'nski in \cite[Theorem 2.2]{kes} determines the shapes of the components in the Auslander-Reiten quiver for any self-injective special biserial algebra of non-polynomial growth.

\begin{theo}[\cite{kes}]\label{tssbatr2} Let $\Lambda$ be a self-injective special biserial algebra of inf\mbox inite-representation type. The following statements are equivalent.
\begin{enumerate}
    \item[$i)$] $\Lambda$ is not an algebra of polynomial growth;
    \item[$ii)$] $\Lambda$ is not a domestic algebra;
    \item[$iii)$] $\leftindex_s\Gamma_{\Lambda}$ consists of a f\mbox inite number of exceptional tubes $\mathbb{Z}A_{\infty}/\langle \tau^r\rangle$, inf\mbox initely many exceptional components of the form $\mathbb{Z}A_{\infty}^{\infty}$, and inf\mbox initely many homogeneous components of the form $\mathbb{Z}A_{\infty}/\langle \tau\rangle$.
\end{enumerate}
\end{theo}

We refer the reader to \cite{kes} and \cite{hpr} for the shapes of the quivers $\mathbb{Z}A_{\infty}/\langle \tau^r\rangle$ and $\mathbb{Z}A_{\infty}^{\infty}$. In the following theorem, the Brauer graphs inducing algebras of polynomial growth are characterized.

\begin{theo}[\cite{bs}]\label{ttrbga} Let $\mathcal{G}$ a Brauer graph and $\Lambda_{\mathcal{G}}$ its associated Brauer graph algebra.
\begin{enumerate}
    \item[$i)$] $\Lambda_{\mathcal{G}}$ is a $1$-domestic algebra if and only if one of the following conditions holds.
    \begin{itemize}
        \item[$(i.a)$] $\mathcal{G}$ is a tree, there are two vertices of $\mathcal{G}$ with multiplicity $2$, and the other vertices have multiplicity $1$;
        \item[$(i.b)$] $\mathcal{G}$ contains a unique cycle, this cycle has odd length, and all the vertices of $\mathcal{G}$ have multiplicity $1$.
    \end{itemize}
    \item[$ii)$] $\Lambda_{\mathcal{G}}$ is a $2$-domestic algebra if and only if $\mathcal{G}$ contains a unique cycle, this cycle has even length, and all the vertices of $\mathcal{G}$ have multiplicity $1$.
    \item[$iii)$] There are no $d$-domestic Brauer graph algebras for $d>2$.
\end{enumerate}
\end{theo}

\subsubsection{Location of string modules}\label{sslsm}

We present some results given by Duf\mbox f\mbox ield in \cite{duff} to describe the \textit{exceptional tubes} of a Brauer graph algebra. Exceptional tubes are the components of the stable Auslander-Reiten quiver containing periodic string modules. We will use these results to study periodic string modules with stable endomorphism ring $\Bbbk$ in Section \ref{sec:main}.

Following \cite{skw}, a module $M$ over a f\mbox inite-dimensional algebra $\Lambda$ is said to be a \textit{periodic module} if $\Omega^n(M)\cong M$, for some positive integer $n$. Here, $\Omega^n(M):=\ker d_{n-1}$ is the \textit{$n$-th syzygy} of $M$, where it is def\mbox ined from a minimal projective resolution of $M$:
$$
\cdots\longrightarrow P_{n+1}\stackrel{d_{n+1}}{\longrightarrow} P_n\longrightarrow\cdots\longrightarrow P_1\stackrel{d_1}{\longrightarrow} P_0 \stackrel{d_0}{\longrightarrow} V\longrightarrow 0
$$
The smallest positive integer $d$ with $\Omega^d (V)\cong V$ is called the \textit{period of} $V$. Now, if $\Lambda$ is a self-injective special biserial algebra, then a string module $M$ is periodic if and only if $M$ corresponds to a vertex in an exceptional tube of $\leftindex_s\Gamma_{\Lambda}$. In particular, it also holds for symmetric special biserial algebras.

\begin{defi}[\cite{duff}]\label{dgwalk} Let $\mathcal{G}=(V,E,\mathfrak{m},\mathfrak{o})$ be a Brauer graph.
\begin{enumerate}
 \item[$i)$] A \emph{Green walk} is a sequence of half edges $((i_j)_{\ast}^{v_j})_{j\geq 0}$, where $v_j\in V$, $(i_j)_{\ast}$ is a half-edge incident to $v_j$ and $\overline{(i_{j+1})_{\ast}^{v_{j+1}}}$ is the successor of $(i_j)_{\ast}^{v_j}$ in the cyclic order of $v_j$. We say that two Green walks $\omega_1$ and $\omega_2$ are the same if $\omega_1$ is a cyclic permutation of $\omega_2$.
 \item[$ii)$] A \emph{double-stepped Green walk} is a subsequence $((i_{2j})_{\ast}^{v_{2j}})_{j\geq 0}$ of a Green walk $((i_j)_{\ast}^{v_j})_{j\geq 0}$.
\end{enumerate}
\end{defi}

\begin{theo}[\cite{duff}]\label{ttgw} Let $\mathcal{G}$ be a Brauer graph and $\Lambda_{\mathcal{G}}$ be its associated Brauer graph algebra. The exceptional tubes in $\leftindex_{s}\Gamma_{\Lambda}$ are in bijection with the double-stepped Green walks of $\mathcal{G}$. Moreover, if $\Delta$ is an exceptional tube corresponding to a Green walk $\omega$, then the rank of $\Delta$ is the length of $\omega$.
\end{theo}

Since the stable Auslander-Reiten quiver of a Brauer graph algebra $\Lambda_{\mathcal{G}}$, with $\mathcal{G}$ its Brauer graph, can be obtained by adding hooks and removing cohooks, we are interested in knowing the modules lying in the boundaries of exceptional tubes. Let $\Lambda_{\mathcal{G}}$ be a Brauer graph algebra.

\begin{defi}[\cite{duff}]\label{dsimuni} We say that a simple $\Lambda_{\mathcal{G}}$-module is a \emph{uniserial simple module} if its projective cover is a uniserial module.
\end{defi}

\begin{theo}[\cite{duff}]\label{tmxu} The modules appearing in the boundary of an exceptional tube of $\leftindex_s\Gamma_{\Lambda_{\mathcal{G}}}$ are the uniserial simple modules and the maximal uniserial submodules of projective indecomposable modules in\hspace{-2.5mm} $\mod\hspace{-2mm} -\hspace{-.5mm} \Lambda_{\mathcal{G}}$.
\end{theo}

The following results are criteria to decide when a simple module and the radical of a projective indecomposable module are in the same component of the stable Auslander-Reiten quiver.

\begin{defi}[\cite{duff}]\label{dexedg} Let $\mathcal{G}=(V,E,\mathfrak{m},\mathfrak{o})$ be a Brauer graph which is not a Brauer tree. An edge $i\in E$ is said to be an \emph{exceptional edge} if there exists a subgraph $T$ of $\mathcal{G}$, called \emph{exceptional subtree of} $\mathcal{G}$, satisf\mbox ing the following conditions.
\begin{enumerate}[$i)$]
    \item The graph $T$ is a tree that contains $i$.
    \item The graph $T$ has a unique vertex $v$, called \emph{connecting vertex}, such that $(\mathcal{G}\backslash T)\cup \{v\}$ is connected.
    \item The graph $T\backslash \{v\}$ shares no vertex with any simple cycle of $\mathcal{G}$. We mean by \emph{simple cycle} to a cycle which is not-crossing at vertices, except at their extreme vertices.
    \item It holds $\mathfrak{m}(w)=1$ for every vertex $w$ of $T\backslash \{v\}$.
\end{enumerate}
\end{defi}

\begin{theo}[\cite{duff}]\label{textub} Assume that $\Lambda_{\mathcal{G}}$ is of inf\mbox inite representation type and let $i$ be an edge of $\mathcal{G}$. Then, the modules $S(i)$ and $\rad(P(i))$ belong to the same exceptional tube in $\leftindex_s\Gamma_{\Lambda}$ if and only if $i$ is an exceptional edge.
\end{theo}

\begin{coro}[\cite{duff}]\label{cexsametub} If $\Lambda_{\mathcal{G}}$ is as in Theorem \ref{textub}, then $S(i)$ and $\rad(P(i))$ are in the same exceptional tube if and only if the half-edges of $i$ are in the same double-stepped Green walk.
\end{coro}

\begin{theo}[\cite{duff}]\label{tstcsrad} Let $i$ and $j$ be (not necessary distinct) non-exceptional edges of $\mathcal{G}$. Then, the modules $S(i)$ and $\rad(P(j))$ are in the same connected component of $\leftindex_s\Gamma_{\Lambda_{\mathcal{G}}}$ if and only if $\Lambda_{\mathcal{G}}$ is $1$-domestic or there is a path 
of even length in $\leftindex_s\Gamma_{\Lambda_{\mathcal{G}}}$ such that the following conditions are satisf\mbox ied.
\begin{enumerate}
   \item[$i)$] The edge $a_k$ is not a loop for every $k\in \{1,\ldots,n\}$.
    \item[$ii)$] For every $k\in \{1,\ldots,n-1\}$, we have $m(v_k)=1$ if $a_k\neq a_{k+1}$ and $m(v_k)=2$ if $a_k=a_{k+1}$.
    \item[$iii)$] For every $k\in \{1,\ldots,n-1\}$, we have that $a_k$ and $a_{k+1}$ are the unique non-exceptional edges incident to the vertex $v_k$.
\end{enumerate}
\end{theo}

We end by pointing out that it is useful to consider a \textit{diagonal} in a component of the stable Auslander-Reiten quiver of a Brauer graph algebra. 

\begin{defi}\label{def:diagonal}
Let $\mathfrak{C}$ be a connected component of $\leftindex_s\Gamma_{\Lambda_{\mathcal{G}}}$. By a \emph{diagonal} in $\mathfrak{C}$ we mean a subgraph $\mathcal{D}$ of $\mathfrak{C}$ of type $A_{\infty}$ such that every path in $\mathcal{D}$ is sectional. 
\end{defi}

Since $\Omega$ is a self-equivalence of $\smod-\Lambda_{\mathcal{G}}$, the image of a diagonal under $\Omega$ is another diagonal. Moreover, if $\mathfrak{C}$ is an exceptional tube or a $\Omega$-stable component of type $\mathbb{Z}A_{\infty}^{\infty}$ for which the image of a diagonal is a non-parallel diagonal, then there exists a diagonal $\mathcal{D}$ of $\mathfrak{C}$ whose vertices are a complete family of representatives of the $\Omega$-orbits of modules in $\mathfrak{C}$. If $\mathfrak{C}$ is an exceptional tube, then $\mathcal{D}$ starts at a module in the boundary of $\mathfrak{C}$.

\subsubsection{Canonical homomorphisms of string modules}

To end this subsection, we recall the main ideas and a result due to H. Krause about canonical homomorphisms between string modules.

\begin{defi}[\cite{kra}]\label{dch}
Let $C$ and $D$ be strings for $\Lambda$. Suppose that $S$ is a substring of both $C$ and $D$ such that the following conditions hold.
\begin{enumerate}
   \item[$i)$] $C\sim ESB$, where $B$ is a substring which is either a trivial path or $B=\alpha B'$ for some arrow $\alpha$, and $E$ is a substring which is a trivial path or $E=E'\beta^{-1}$ for some arrow $\beta$;
    \item[$ii)$] $D\sim FSG$, where $G$ is a substring which is either a trivial path or $G=\gamma^{-1}G'$ for some arrow $\gamma$, and $F$ is a substring which is either a trivial path or $F=F'\delta$ for some arrow $\delta$.
\end{enumerate}
In this situation, we have an epimorphism $M[C]\twoheadrightarrow M[S]$ and a monomorphism $M[S]\xhookrightarrow{} M[D]$. We def\mbox ine the \emph{canonical homomorphism} induced by the factorizations $ESB$ and $FSG$ as the composition $M[C]\twoheadrightarrow M[S]\xhookrightarrow{} M[D]$.
\end{defi}

The following result is due to H. Krause in \cite{kra}.

\begin{theo}[\cite{kra}]\label{tcnn} Every homomorphism of string modules is a linear combination of canonical homomorphisms.
\end{theo}

Sometimes, we will work with canonical homomorphisms between a string module and a projective indecomposable module. For a string module $M[C]$ and a projective indecomposable module $P$, by a canonical homomorphism from $M[C]$ to $P$ we mean a composition $\iota f$, where $f:M[C]\rightarrow \rad(P)$ is a canonical homomorphism and $\iota:\rad(P)\xhookrightarrow{} P$ is the natural inclusion. By a canonical homomorphism from $P$ to $M[C]$ we mean a composition $g\varphi$, where $g:M[C]\rightarrow P$ is a canonical homomorphism and $\varphi:P\twoheadrightarrow P/\soc(P)$ is the canonical projection.

\section{Periodic string modules and its universal deformation rings}\label{sec:main}

In this section, we determine the universal deformation rings of periodic string modules with stable endomorphism ring isomorphic to $\Bbbk$ over generalized Brauer tree algebras of non-polynomial growth.

By Theorem \ref{tssbatr2} and Theorem \ref{ttrbga}, we can assume that all the Brauer graphs considered in this paper are generalized Brauer trees with at least two vertices having multiplicities greater than $1$ and we can exclude the generalized Brauer trees having exactly two vertices with multiplicity $2$ and all the other vertices having multiplicity $1$.

\subsection{A derived equivalence involving a star}\label{sedis}

This subsection introduces a specif\mbox ic generalized Brauer graph, a star, which will be shown to be central to this work.

\begin{figure}[h!]
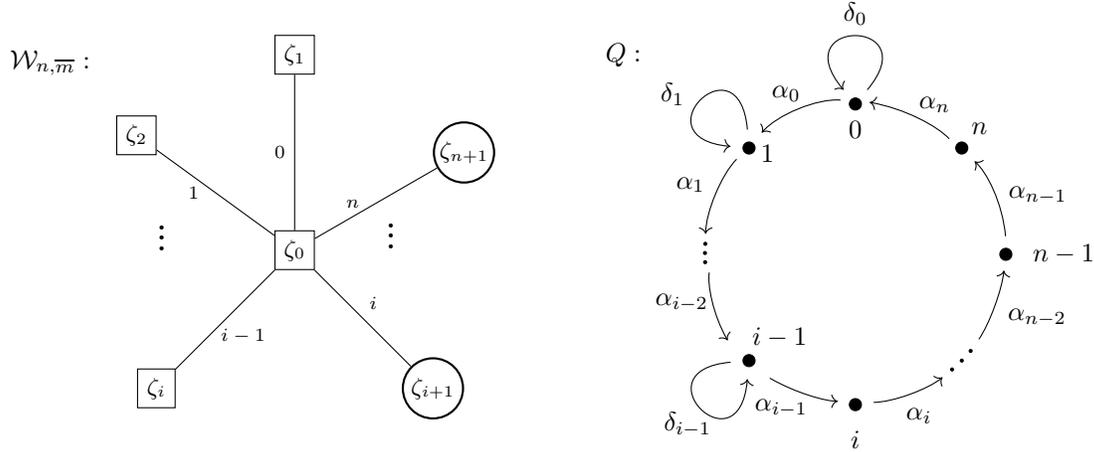

\mypicture\hspace{10mm}
\mypictureb[baseline=-22mm]
\caption{The star $\mathcal{W}_{n,\overline{m}}$ and its associated quiver $Q$}\label{fig1}
\end{figure}

Let $n$ and $i$ be non-negative integers, with $0\leq i\leq n+1$. Given a vector $\overline{m}=(m_0,\ldots,m_i)$ of $i+1$ integers such that $2\leq m_0\leq \ldots \leq m_i$, def\mbox ine the \textit{star} $\mathcal{W}_{n,\overline{m}}$ to be the Brauer graph shown in Figure \ref{fig1}, where we draw circles for vertices of multiplicity $1$ and squares for vertices of multiplicity greater than $1$. According to this, the multiplicity of $\zeta_j$ is $m_j$ for all $j\in \{0,\ldots,i\}$ and all the other vertices have multiplicity $1$. We assume that the edges of $\mathcal{W}_{n,\overline{m}}$ are in $\mathbb{Z}/(n+1)\mathbb{Z}$.

We denote by $\Lambda(n,\overline{m})$ the Brauer graph algebra associated to $\mathcal{W}_{n,\overline{m}}$. Then, $\Lambda(n,\overline{m})=\Bbbk Q/I$, where $Q$ is showed in Figure \ref{fig1} and the ideal $I$ is def\mbox ined by 
$$
I=\langle \delta_j^{m_{j+1}}-A_j^{m_0},\alpha_{j-1} \delta_j,\delta_j \alpha_j,A_k^{m_0}\alpha_k\mid
j\in\{0,\ldots,i-1\} \mbox{ and } k\in\{i,\ldots,n-1\} \rangle,
$$ 
with $A_l:=A_{\zeta_0,\hat{l}}=\alpha_l \cdots \alpha_{l-1}$ for all $l\in \mathbb{Z}/(n+1)\mathbb{Z}$.

The stable Auslander-Reiten quiver of $\Lambda(n,\overline{m})$ coincides with the Auslander-Reiten quiver of $\Lambda_s:=\Lambda(n,\overline{m})/\soc(\Lambda(n,\overline{m}))\cong \Bbbk Q/I_s$, where
$$I_s=\langle \delta_j^{m_{j+1}},A_l^{m_0},\alpha_{j-1} \delta_j,\delta_j \alpha_j:j\in\{0,\ldots,i-1\},l\in\{0,\ldots,n\}\rangle.$$

In the following def\mbox inition we introduce a class of strings of $\Lambda_s$ that will appear along this section.

\begin{defi}\label{dsstr} 
For distinct vertices $l$ and $p$ of the quiver $Q$ (as in Figure \ref{fig1}), we def\mbox ine
$$x_{l,l}=e_l,$$
$$x_{l,p}=\alpha_l \cdots \alpha_{p-1}.$$
For $l,p\in \{0,\ldots,i-1\}$, with $l<p$, we def\mbox ine
$$y_l=(\delta_l^{-(m_{l+1}-1)}\alpha_l) \cdots (\delta_{i-1}^{-(m_i-1)}\alpha_{i-1}).$$
\end{defi}

\begin{lemm}\label{lbdes} Let $\Lambda_{\mathcal{G}}$ be a Brauer graph algebra whose underlying Brauer graph is a generalized Brauer tree $\mathcal{G}$. Then, $\Lambda_{\mathcal{G}}$ is derived equivalent to a Brauer graph algebra whose underlying Brauer graph is a star with the same multiset of multiplicities of $\mathcal{G}$. In particular, there exist a bijective correspondence between exceptional tubes that preserves ranks, modules with stable endomorphism ring isomorphic to $\Bbbk$ and universal deformation rings.
\end{lemm}

\begin{proof}
Set $\mathcal{G}=(V,E,\mathfrak{m},\mathfrak{o})$. Then $|V|\geq 2$. We put $n=|V|-2$ and $i=|V^{\ast}|-1$, where $V^{\ast}:=\{v\in V: \mathfrak{m}(v)>1\}$. Write $V^{\ast}=\{v_0,\ldots,v_i\}$. Without loss of generality, assume that $\mathfrak{m}(v_0)\leq \ldots \leq \mathfrak{m}(v_i)$. We have that $|E|=|V|-1=n+1$ because $(V,E)$ is a tree. Thus, there exists a bijective correspondence between the vertices, the edges, and the multiset of multiplicities of $\mathcal{G}$ and $\mathcal{W}_{n,\overline{m}}$, where $\overline{m}=(\mathfrak{m}(v_0),\ldots,\mathfrak{m}(v_i))$. Moreover, the fact that $(V,E)$ contains no cycles implies that both generalized Brauer graphs, $\mathcal{G}$ and $\mathcal{W}_{n,\overline{m}}$, determine a unique face whose perimeter is $2\cdot(n+1)=2n+2$. Since every tree is bipartite, we have that $(V,E)$ and $\mathcal{W}_{n,\overline{m}}$ are bipartite.

In conclusion, by Opper-Zvonareva's theorem (see \cite[Theorem A]{oz}), we have that $\Lambda_{\mathcal{G}}$ and $\Lambda(n,\overline{m})$ are derived equivalent.
\end{proof}

\subsection{Exceptional tubes}\label{sub:exctub}

Let $\Lambda_{\mathcal{G}}$ be a Brauer graph algebra of non-polynomial growth whose underlying Brauer graph is a generalized Brauer tree $\mathcal{G}=(V,E,\mathfrak{m},\mathfrak{o})$. In this subsection, we determine the exceptional tubes appearing in the stable Auslander-Reiten quiver of $\mathcal{G}$. Recall that, by Theorem \ref{tssbatr2}, the components of $\leftindex_s\Gamma_{\Lambda_{\mathcal{G}}}$ are exceptional tubes, components of type $\mathbb{Z}A_{\infty}^{\infty}$ and homogeneous tubes.

In the following lemma, we determine the number of exceptional tubes and their ranks. We will use Theorem \ref{ttgw} which gives a criterion to determine these numbers depending on the concept of \textit{double-stepped Green walk} described in Def\mbox inition \ref{dgwalk} $ii)$. 

\begin{lemm}\label{lcgbt} The stable Auslander-Reiten quiver of $\Lambda_{\mathcal{G}}$ is a disjoint union of $2$ exceptional tubes of rank $|E|$, inf\mbox initely many (exceptional) components of the form $\mathbb{Z} A_{\infty}^{\infty}$, and inf\mbox initely many homogeneous tubes.
\end{lemm}

\begin{proof} By Theorem \ref{tssbatr2}, we know that $\leftindex_s \Gamma_{\Lambda_{\mathcal{G}}}$ is the disjoint union of a f\mbox inite number of exceptional tubes, inf\mbox initely many homogeneous tubes and inf\mbox initely many exceptional components of the form $\mathbb{Z}A_{\infty}^{\infty}$. By Theorem \ref{ttgw}, the number of exceptional tubes in $\leftindex_s \Gamma_{\Lambda_{\mathcal{G}}}$ equals the number of double-stepped Green walks of $\mathcal{G}=(V,E,\mathfrak{m},\mathfrak{o})$. Thus, it is enough to verify that $\mathcal{G}$ just admits two double-stepped Green walks.

With this goal in mind, f\mbox irstly, we will prove by induction on $|E|$ that $\mathcal{G}$ only has one Green walk (see Def\mbox inition \ref{dgwalk} $i)$). The case in which $\mathcal{G}$ has only an edge is clear. Suppose that $|E|>1$ and that every generalized Brauer tree with $|E|-1$ edges just admits one Green walk. Let $v$ be a vertex of valency $1$ in $\mathcal{G}$, $y(v)$ the unique edge adjacent to $v$, and $v'$ the other vertex in which $y(v)$ incides. Moreover, we denote by $x_0$ the successor of $y(v)$ in the cyclic order given by $v'$. Since $|E|>1$ and $\mathcal{G}$ does not contain cycles, we have $x_0\neq y(v)$. Let $\mathcal{G}\backslash v$ be the Brauer graph obtained from $\mathcal{G}$ by removing the vertex $v$, the edge $y(v)$, and with multiplicity and orientation functions def\mbox ined as the corresponding restrictions of $\mathfrak{m}$ and $\mathfrak{o}$. Observe that if the cyclic order at $v'$ in $\mathcal{G}$ is given by $y(v)<x_0<...<x_s<y(s)$, then the cyclic order at $v'$ in $\mathcal{G}\backslash v$ is given by $x_0<...<x_s<x_0$. Since $\mathcal{G}\backslash v$ is a generalized Brauer graph with $|E|-1$ edges, from the induction hypothesis we obtain that $\mathcal{G}\backslash v$ only admits one Green walk $g$. Then, $g$ contains to every half-edge of $\mathcal{G}\backslash v$ and, in particular, we can write $g=(x_0^{v_0},z_1^{v_1},\ldots,z_h^{v_h},x_s^{v'})$ where $v_0$ is the vertex in $V$ such that $v_0\neq v'$ and $x_0$ is incident to $v_0$. Since $x_0=z_1$ is the successor of $y(v)$ and $x_s=z_h$ is the predecessor of $y(s)$, both in the cyclic order of $v'$, it follows that none of the other edges $z_2,\ldots,z_{h-1}$ is the successor or the predecessor of $y(v)$. Thus, the f\mbox inite sequences $(z_1^{v_1},\ldots,z_h^{v_h})$ and $(x_s^{v'},y(v)^v,y(v)^{v'},x_0^{v_0})$ are subsequences of Green walks in $\mathcal{G}$. These two Green walks have to be the same due to the fact that they share common half-edges, and this is given by $g'=(y(v)^v,y(s)^{v'},x_0^{v_0},z_1^{v_1},\ldots,z_h^{v_h},x_s^{v'},y(v)^v)$. Indeed, it is the unique Green walk in $G$, because it contains all the half-edges of $\mathcal{G}$: all the half-edges of $\mathcal{G}\backslash v$ joint with $y(v)^v$ and $y(v)^{v'}$, which concludes the proof of the existence of the unique Green walk.

Now, for simplicity, we denote the only Green walk of $\mathcal{G}$ by $g'=(g_1^{u_1},\ldots,g_r^{u_r})$. Since there is an even number of half-edges, $r$ is even, and we have two double-stepped Green walks, namely, $(g_1^{u_1},g_3^{u_3},\ldots,g_{r-1}^{u_{r-1}})$ and $(g_2^{u_2},g_4^{u_4},\ldots,g_r^{u_r})$. Each double-stepped Green walk has length $r/2=|E|$, because $r$ coincides with the number of half-edges of $\mathcal{G}$.
\end{proof}

Now, our current focus is on understanding the action of the syzygy operator on exceptional tubes. In the following lemma we state that the two tubes of $\leftindex_s\Gamma_{\Lambda_{\mathcal{G}}}$ are in the same $\Omega$-orbit.

\begin{lemm}\label{ltubsyz} If $\mathfrak{C}_1$, $\mathfrak{C}_2$ are the exceptional tubes of $\leftindex_s\Gamma_{\Lambda_{\mathcal{G}}}$ then, $\Omega(\mathfrak{C}_1)=\mathfrak{C}_2$ and $\Omega(\mathfrak{C}_2)=\mathfrak{C}_1$.
\end{lemm}

\begin{proof}
Let $\mathcal{G}$ be the generalized Brauer tree associated to $\Lambda_{\mathcal{G}}$ and $\mathcal{W}_{n,\overline{m}}$ be the star for which the algebras $\Lambda_{\mathcal{G}}$ and $\Lambda(n,\overline{m})$ are derived equivalent. Then, by Theorem \ref{tderst}, there is an equivalence $F:\smod-\Lambda_{\mathcal{G}} \rightarrow \smod-\Lambda(n,\overline{m})$. In particular, $F$ induces a bijection between the exceptional tubes of $\leftindex_s\Gamma_{\Lambda_{\mathcal{G}}}$ and $\leftindex_s\Gamma_{\Lambda(n,\overline{m})}$ preserving ranks. Let $\mathfrak{C}_1$ and $\mathfrak{C}_2$ be the exceptional tubes of $\leftindex_s\Gamma_{\Lambda_{\mathcal{G}}}$, and $\mathfrak{C}_3=F(\mathfrak{C}_1)$ and $\mathfrak{C}_4=F(\mathfrak{C}_2)$ the exceptional tubes of $\leftindex_s\Gamma_{\Lambda(n,\overline{m})}$.

First, we will see that the statement holds for $\mathcal{W}_{n,\overline{m}}$ and after this we apply the equivalence $F$ to get the desired result for any generalized Brauer tree.

Let $j$ be an exceptional edge in $\mathcal{W}_{n,\overline{m}}$. Without loss of generality, suppose that $S(j)$ belongs to $\mathfrak{C}_3$. By the construction of the two double-stepped Green walks for generalized Brauer trees in the proof of Lemma \ref{lcgbt}, the half-edges of $j$ are in dif\mbox ferent double-stepped Green walks. Hence, by Corollary \ref{cexsametub}, the modules $S(j)$ and $\rad(P(j))$ are in distinct connected components. Since $\Omega(S(j))=\rad(P(j))$, it follows that the components containing to $S(i)$ and $\Omega(S(i))$ are dif\mbox ferent. Thus, it holds $\Omega(\mathfrak{C}_3)=\mathfrak{C}_4$ and $\Omega(\mathfrak{C}_4)=\mathfrak{C}_3$.

Now, suppose that $\Omega(\mathfrak{C}_1)=\mathfrak{C}_1$. Then, there is a module $M$ in $\mathfrak{C}_1$ and an undirected path connecting $M$ and $\Omega(M)$ in $\mathfrak{C}_1$. Applying the equivalence $F$, we get an undirected path between $F(M)$ and $\Omega(F(M))$ in $\leftindex_s\Gamma_{\Lambda_{\mathcal{G}}}$, which is a contradiction with the fact that $\mathfrak{C}_3\neq \Omega(\mathfrak{C}_3)$. So, it holds $\Omega(\mathfrak{C}_1)=\mathfrak{C}_2$. Since $\Omega$ is a self-equivalence of the stable category, it also follows that $\Omega(\mathfrak{C}_2)=\mathfrak{C}_1$.
\end{proof}

\begin{remark}\label{rexedst} \rm The unique simple modules in exceptional tubes for the star $\mathcal{W}_{n,\overline{m}}$ are the simple modules corresponding to the edges incident to a truncated vertex. To see this, we verify that the exceptional edges in $\mathcal{W}_{n,\overline{m}}$ are the edges incident to a truncated vertex. In fact, note that all the vertices in which the edges $0$, ..., $i-1$ are incident have multiplicities greater than $1$, and then they do not belong to any exceptional subtree of $\mathcal{W}_{n,\overline{m}}$. On the other hand, for every edge $j\in \{i,\ldots,n\}$, the subgraph of $\mathcal{W}_{n,\overline{m}}$ formed only by $j$ is an exceptional subtree of $\mathcal{W}_{n,\overline{m}}$ with connecting vertex $\zeta_0$. Therefore, all the exceptional edges of $\mathcal{W}_{n,\overline{m}}$ are $i$, ..., $n$. By Theorem \ref{textub}, these are the edges corresponding to simple modules in the exceptional tubes of $\leftindex_s\Gamma_{\Lambda(n,\overline{m})}$.
\end{remark}

\subsection{Universal deformation rings of periodic modules}

Let $\Lambda_{\mathcal{G}}$ be a Brauer graph algebra whose underlying Brauer graph $\mathcal{G}=(V,E,\mathfrak{m},\mathfrak{o})$ is a generalized Brauer tree. To compute universal deformation rings of modules over $\Lambda_{\mathcal{G}}$, we use the derived equivalence provided in Lemma \ref{lbdes} to reduce the calculation to the case for the star $\mathcal{W}_{n,\overline{m}}$, with $n=|V|-2$ and $\overline{m}=(m(v_0),\ldots,m(v_i))$, where $i+1$ is the number of non-truncated vertices in $V$. Recall $\Lambda(n,\overline{m})$ denotes the Brauer graph algebra associated to $\mathcal{W}_{n,\overline{m}}$.

We will use the next def\mbox inition to describe the periodic string modules with stable endomorphism ring isomorphic to $\Bbbk$.

\begin{defi}\label{ddm} Consider $\Lambda$ a self-injective special biserial algebra. Let $\mathfrak{C}$ be an exceptional tube in the stable Auslander-Reiten quiver of $\Lambda$ and $M$ be a string module in $\mathfrak{C}$. We def\mbox ine the \emph{distance} from $M$ to the boundary of $\mathfrak{C}$ as the non-negative integer $d_M$ for which there exists a sectional path of the form
$$\xymatrix{M_{d_M} \ar[r] & M_{d_M-1} \ar[r] & \cdots \ar[r] & M_0},$$
where $M_0=M$ and $M_t$ is a module in the boundary of $\mathfrak{C}$. If $M$ belongs to the boundary of $\mathfrak{C}$, then we def\mbox ine $d_M:=0$.
\end{defi}

We consider two cases based on the existence of vertices with multiplicity 1. While the results and proofs are similar in both cases, we divide them for technical reasons related to the shapes of projective modules and, consequently, the shapes of the strings of $\Lambda(n,\overline{m})/\soc(\Lambda(n,\overline{m}))$.

\subsubsection{First case of the classif\mbox ication}\label{subsec:one}

In this subsection, we assume that at least one vertex in $\mathcal{W}_{n,\overline{m}}$ has multiplicity $1$, that is $i\neq n+1$. The radical series of the projective indecomposable modules over $\Lambda(n,\overline{m})$ are exhibited in Figure \ref{fig2}.
\begin{figure}[h!]
\centering
\begin{subfigure}[b]{0.6\textwidth}
\centering
\def\dist{1.5}
$P(j):$\ \ \tabbedCenterstack{& $S(j)$ &&& \\
$S(j)$ && $S(j+1)$ && \\
 &&&& \\
\vdots && \vdots  && \\
 && $S(j)$ && \\
$S(j)$ &&  && \\
 && $S(j+1)$ && \mbox{\ \ , for $j\in\{0,\ldots,i-1\}$;}\\
&&&& \\
\vdots && \vdots && \\
&&&& \\
$S(j)$ && $S(j-1)$ && \\
& $S(j)$ &&& }
\end{subfigure}
\begin{subfigure}[b]{0.4\textwidth}
\centering
\def\dist{1.5}
$P(j):$\ \ \tabbedCenterstack{& $S(j)$\\
& \\
& $S(j)$\\
& \\
& \vdots\\
& \\
& $S(j)$}
\mbox{\ \ , for } $j\in\{i,\ldots,n\}$.
\end{subfigure}
\caption{Radical series of projective indecomposable modules in $\mathcal{W}_{n,\overline{m}}$.}\label{fig2}
\end{figure}

By Theorem \ref{tmxu}, the modules appearing in the boundaries of exceptional tubes are the uniserial simple modules and the maximal uniserial submodules of projective indecomposable modules. These modules are $M[x_{l+1,l} \cdot A_l^{m_0-1}]$, $M[\delta_j^{m_{j+1}-1}]$ and $S(k)$, with $l\in \{0,\ldots, n\}$, $j\in \{0,\ldots, i-1\}$, and $k\in \{i,\ldots, n\}$. See Def\mbox inition \ref{dsstr} to recall the def\mbox inition of the strings $x_{l,p}$.

\begin{lemm}\label{lsmmp} Let $\mathcal{W}_{n,\overline{m}}$ be a star, with $\overline{m}=(m_0,\ldots,m_i)$. Suppose that at least one vertex of $\mathcal{W}_{n,\overline{m}}$ has multiplicity $1$. Let $M$ be a string module in an excepcional tube $\mathfrak{C}$ of $\leftindex_s \Gamma_{\Lambda(n,\overline{m})}$. Then, $d_M\leq n$ if and only if $M$ belongs to the $\Omega$-orbit of some string module of the form $M[C]$, where $C$ is one of the following strings:
\begin{enumerate}
    \item[$i)$] $x_{k,n}$, with $k\in \{i,i+1,\ldots,n\}$;
    \item[$ii)$] $y_j \cdot x_{i,n}$, with $j\in \{0,\ldots,i-1\}$.
\end{enumerate}
Moreover, $d_M=n$ if and only if $M$ is in the $\Omega$-orbit of the string module $M[y_0\cdot x_{i,n}]$.
\end{lemm}
\begin{proof}
By Lemma \ref{lcgbt} (and Lemma \ref{lbdes}), there exist just two exceptional tubes $\mathfrak{C}_1$ and $\mathfrak{C}_2$ in $\leftindex_s \Gamma_{\Lambda(n,\overline{m})}$, each one being the image of the other via the self-equivalence $\Omega$. More precisely, this self-equivalence induces an isomorphism of these two tubes seen as valued translation quivers. In particular, the irreducible homomorphisms from $M$ (resp. to $M$) are in bijective correspondence with the irreducible homomorphisms from $\Omega(M)$ (resp. to $\Omega(M)$) for every non-projective indecomposable module $M$. Thus, if $D_1$ and $D_2$ denote the subgraphs of $\mathfrak{C}_1$ and $\mathfrak{C}_2$, respectively, consisting of the modules $M$ for which $d_M\leq n$ and the irreducible homomorphisms between them, then $D_1=\Omega(D_2)$ and $D_2=\Omega(D_1)$. Hence, it is enough to consider only one exceptional tube.

Let $\mathfrak{C}$ be the exceptional tube containing to the uniserial simple module $S(n)$ and let $\mathcal{D}$ be the diagonal of $\mathfrak{C}$ starting at $S(n)$ (see Def\mbox inition \ref{def:diagonal}). Note that all the modules in $\mathfrak{C}$ are in the $\tau$-orbits of $\mathcal{D}$, and hence in the $\Omega$-orbits of $\mathcal{D}$ because $\tau$ and $\Omega^2$ are naturally isomorphic functors. Thus, it only remains to verify that the modules $M$ of $\mathcal{D}$ for which $d_M\leq n$ are those given in the statement. By inspection of the construction of $\mathfrak{C}$ by adding and removing hooks and cohooks, we obtain that the f\mbox irst $n+1$ terms of the diagonal of $\mathfrak{C}$ starting at $S(n)$ are the following
\[
\xymatrix{
M_n \ar[r] & M_{n-1} \ar[r] & \cdots \ar[r] & M_i \ar[r] & M_{i-1} \ar[r] & \cdots \ar[r] & M_1 \ar[r] & M_0\ ,
}
\]
with $M_n:=S(n)=M[x_{n,n}]$ and $M_{l-1}=(M_l)_h$ for $l\in \{1\ldots,n\}$, which gives
$$M_k:=M[x_{k,n}] \mbox{ for } k\in \{i,\ldots,n-1\} \mbox{ and }$$
$$M_j:= M[y_j\cdot x_{i,n}] \mbox{ for } j\in \{0,\ldots,i-1\}.$$
\end{proof}

\begin{remark}\label{rd} 
\rm
The complete diagonal of $\mathfrak{C}$ starting at $S(n)$ is given by
\[
\xymatrix{
M_n \ar[r] & \cdots \ar[r] & M_0 \ar[r] & N_0 \ar[r] & N_1 \ar[r] & \cdots
}
\]
where the modules $M_l$, with $l\in \{0,\ldots,n\}$, are as in the proof of Lemma \ref{lsmmp}, $N_0=(M_0)_h$ and $N_{p+1}=(N_p)_h$ for $p\in \{0,1,\ldots\}$. More precisely, if $p=q(n+1)+r$, with $q,r\in \mathbb{Z}$ and $0\leq r\leq n$, we put $\rho:=(\alpha_n y_0 x_{i,n})^{q+1}$ and we def\mbox ine the string $C_p$ as follows:
\begin{itemize}
    \item[$i)$] If $r=0$, then $C_p=\rho$.
    \item[$ii)$] If $0<r\leq n-i$, then $C_p=x_{n-r,n} \rho$.
    \item[$iii)$] If $n-i<r\leq n$, then $C_p=y_{n-r} x_{i,n} \rho$.
\end{itemize}
Thus, $N_p=M[C_p]$.
\end{remark}

\begin{theo}\label{tmek} Suppose that at least one vertex of $\mathcal{W}_{n,\overline{m}}$ has multiplicity $1$. Let $M$ be a periodic string module over $\Lambda(n,\overline{m})$. The stable endomorphism ring of $M$ is isomorphic to $\Bbbk$ if and only if $d_M\leq n$.
\end{theo}
\begin{proof}
As discussed in Subsection \ref{sub:exctub}, periodic string modules coincide with the modules in the exceptional tubes of $\leftindex_s \Gamma_{\Lambda(n,\overline{m})}$. Additionally, we can use the string modules $M_j$ and $N_p$ as described in Remark \ref{rd}, since $\send_{\Lambda(n,\overline{m})}(M)\cong \send_{\Lambda(n,\overline{m})}(\Omega(M))$ for all non-projective indecomposable module $M$ on $\Lambda(n,\overline{m})$. Adopting the notation used in the proof of Lemma \ref{lsmmp} and Remark \ref{rd}, we will now prove that $\send_{\Lambda(n,\overline{m})}(M_l)\cong \Bbbk$ for every $l\in \{0,\ldots,n\}$.

For the string module $M_k=M[x_{k,n}]$, with $k\in \{i,\ldots,n\}$, the unique substring $S$ of $x_{k,n}$ satisfying the Def\mbox inition \ref{dch} is itself, and the unique canonical endomorphism of $M_k$ is the identity. This implies that $\dim_{\Bbbk}(\send_{\Lambda(n,\overline{m})}(M_k))=1$.

For the string modules $M_j=M[y_j x_{i,n}]$, with $j\in \{0,\ldots,i-1\}$, the substrings of $y_j x_{i,n}$ satisfying Def\mbox inition \ref{dch} are $y_j x_{i,n}$ and $\delta_l^t$, where $l\in \{j,\ldots,i-1\}$ and $0\leq t\leq m_{l+1}-2$. Recall that $m_{l+1}$ is the multiplicity of $\zeta_{l+1}$. Each substring $\delta_l^t$ uniquely determines a canonical endomorphism of  $M_j$. Let $f$ be the canonical endomorphism of $M_j$ induced by $\delta_l^t$. We def\mbox ine the strings $D$ and $E$ as follows.
\begin{description}
\item[\textbf{Case 1}] If either $j\leq l=i-1$ and $t<m_i-2$ or $j=l=i-1$ and $t=m_i-2$, then $D:=\delta_{i-1}^{t+1}$ and $E:=\delta_{i-1}^{-(m_i-1)}\alpha_{i-1}\cdots \alpha_{n-1}$.\\
\item[\textbf{Case 2}] If $j<l=i-1$ and $t=m_i-2$, then we set $D:=\alpha_{i-2} \delta_{i-1}^{-(m_i-1)}$ and $E:=\delta_{i-1}^{-(m_i-1)}\alpha_{i-1}\cdots \alpha_{n-1}$.\\
\item[\textbf{Case 3}] If either $j<l<i-1$ and $t<m_{l+1}-2$, or $j=l<i-1$, then $D:=\delta_l^{t+1}$ and $E:=\delta_l^{-(m_{l+1}-1)}\alpha_l$.\\
\item[\textbf{Case 4}] If $j<l<i-1$ and $t=m_{l+1}-2$, then $D:=\alpha_{l-1}\delta_l^{-(m_{l+1}-1)}$ and $E:=\delta_l^{-(m_{l+1}-1)}\alpha_l$.
\end{description}
For substrings $D$ and $E$ of $C$, we have that $M[D]$ is a quotient module of $M_j$ and a submodule of $P(l)$, while $M[E]$ is a quotient module of $P(l)$ and a submodule of $M_j$. The canonical homomorphisms induced by $M[D]$ and $M[E]$, respectively, are $\varphi \in \Hom_{\Lambda(n,\overline{m})}(M_j,P(l))$ and $\psi \in \Hom_{\Lambda(n,\overline{m})}(P(l),M_j)$. Consequently, $f=\psi \varphi$. This implies that the only canonical endomorphism of $M_j$ that does not factor through any projective module is the identity. By Theorem \ref{tcnn}, it follows that $\dim_{\Bbbk}(\send_{\Lambda(n,\overline{m})}(M_j))=1$.

To complete the proof, we show that the remaining modules on the diagonal of $\mathfrak{C}$ starting at $S(n)$ do not have stable endomorphism rings isomorphic to $\Bbbk$. These modules are the $N_p$ def\mbox ined in Remark \ref{rd}. Using the notations in Remark \ref{rd}, we construct an endomorphism that does not factor through any projective module. We def\mbox ine substrings $F$ and $G$ of $C_p$ to obtain a canonical endomorphism $g:N_p\twoheadrightarrow M[F]\xhookrightarrow{} N_p$. . We now consider the following cases.
\begin{description}
\item[\textbf{Case 1'}] If $r=0$, then we set $F:=e_n$ and $G:=\rho$.\\
\item[\textbf{Case 2'}] If $0<r\leq n-i$, then we set $F:=x_{n-r,n}$ and $G:=x_{n-r,n}(\alpha_n y_0 x_{i,n})^q \alpha_n y_0x_{i,n-r}$.\\
\item[\textbf{Case 3'}] If $n-i<r\leq n$, then we set $F:=\alpha_{i-1}x_{i,n}$ and $G:=y_{n-r} x_{i,n}(\alpha_n y_0 x_{i,n})^q \alpha_n y_{0,i-1}\delta_{i-1}^{m_i-1}$.
\end{description}
For the f\mbox irst two cases, $C_p=F \rho$, while for the third case, $C_p=y_{n-r,i-1}\delta_{i-1}^{m_i-1}F \rho$. In all cases, $C_p=G F$. These factorizations induce the desired endomorphism $g$. To show that $g$ does not factor through any projective module, suppose there exists a projective module $P$ and homomorphisms $\varphi \in \Hom_{\Lambda(n,\overline{m})}(N_p,P)$ and $\sigma \in \Hom_{\Lambda(n,\overline{m})}(P,N_p)$ such that $g=\sigma \varphi$. The only canonical homomorphisms from $N_p$ to an indecomposable direct summand of $P$ that preserve the quotient module $M[F]$ of $N_p$ are $\beta$ and $\gamma$ def\mbox ined as follows: If $P$ has a direct summand of the form $P(n)$, $\beta$ is the canonical homomorphism $N_p\twoheadrightarrow M[F] \xhookrightarrow{} P(n)$; otherwise, $\gamma$ is any canonical homomorphism $N_p\twoheadrightarrow M[A] \xhookrightarrow{} P(0)$, where $A=F\alpha_nA'$ for some string $A'$. Due to the lengths of the radical series of $P(n)$ and $P(0)$, we obtain that $\sigma \iota_n \beta=0$ and $\sigma \iota_0 \gamma=0$, where $\iota_n$ and $\iota_0$ are the natural inclusions $P(n)\xhookrightarrow{} P$ and $P(0)\xhookrightarrow{} P$, respectively. Since $\varphi$ is a linear combination of canonical homomorphisms, $M[F]$ maps to $\ker \sigma$ under $\varphi$, contradicting the fact that $g=\sigma \varphi$. Therefore, $g$ does not factor through a projective module, and by Theorem \ref{tcnn}, $\dim_{\Bbbk} \send_{\Lambda(n,\overline{m})}(M[N_p])>1$. This completes the proof.
\end{proof}

\begin{theo}\label{text} Suppose that at least one vertex of $\mathcal{W}_{n,\overline{m}}$ has multiplicity $1$. Let $M$ be a periodic string $\Lambda(n,\overline{m})$-module such that $\send_{\Lambda(n,\overline{m})}(M)\cong \Bbbk$. Then,
\vspace{-3.5mm}

\begin{equation*}
\dim_{\Bbbk}\Ext_{\Lambda(n,\overline{m})}^1(M,M)= \begin{cases}
0 & \text{if $d_M<n$,}\\
1 & \text{if $d_M=n$.}
\end{cases}
\end{equation*}
\end{theo}

\begin{proof}
Since $\Omega$ is a self-equivalence on $\smod-\Lambda(n,\overline{m})$ we have that, for each indecomposable 
non-projective module $M$,
$$
\Ext_{\Lambda(n,\overline{m})}^1(M,M)\cong \shom_{\Lambda(n,\overline{m})}(\Omega(M),M)\cong \shom_{\Lambda(n,\overline{m})}(\Omega^2(M),\Omega(M)) \cong \Ext_{\Lambda(n,\overline{m})}^1(\Omega(M),\Omega(M)).
$$
In consequence, we can focus on the modules $M_l$ in the proof of Lemma \ref{lsmmp}, where $l\in \{0,\ldots,n\}$. Thus, it suf\mbox f\mbox ices to show that $\Ext_{\Lambda(n,\overline{m})}(M_l,M_l)=0$ for $l\neq 0$ and $\dim_{\Bbbk}\Ext_{\Lambda(n,\overline{m})}(M_0,M_0)=1$. We consider the following three cases.
\begin{description}
\item[\textbf{Case 1}] $l\in \{i,\ldots,n\}$, or $i\neq 1$ and $l=i-1$.\\
For $l\in \{i,\ldots,n\}$, $M_l=M[x_{l,n}]$, and for $i\neq 1$ and $l=i-1$, then $M_{i-1}=M[y_{i-1}\cdot x_{i,n}]$. In both cases, taking the kernel of the projective cover $P(l)\twoheadrightarrow M_l$, we obtain $\Omega(M_l)=M[x_{0,l}A_l^{m_0-1}]$. Consequently, there exists no string $C$ inducing a canonical homomorphism of the form $\Omega(M_l) \twoheadrightarrow M[C] \xhookrightarrow{} M_l$. By Theorem \ref{tcnn}, this implies that $\Hom_{\Lambda(n,\overline{m})}(\Omega(M_l),M_l)=0$ and thus
$$\dim_{\Bbbk} \Ext_{\Lambda(n,\overline{m})}^1(M_l,M_l)=\dim_{\Bbbk} \shom_{\Lambda(n,\overline{m})}(\Omega(M_l),M_l)=0.$$

\item[\textbf{Case 2}] $i\notin \{1,2\}$ and $l\in \{1,\ldots,i-2\}$.\\
We have that $M_l=M[y_l\cdot x_{i,n}]$. Then, using the projective cover
$$P(l)\oplus P(l+1) \oplus \cdots \oplus P(i-1)\twoheadrightarrow M_l,$$
we obtain $\Omega(M_l)=M[S]$, where
$$
S=x_{0,i-1}A_{i-1}^{m_0-1} \cdot x_{i-1,i-2}A_{i-2}^{m_0-1}\cdot \ldots \cdot x_{l+1,l}A_l^{m_0-1}.
$$
We have that the unique strings $C$ inducing a canonical endomorphism $\Omega(M_l) \twoheadrightarrow M[C] \xhookrightarrow{} M_l$ are $e_{l+1}$, $e_{l+2}$, ..., $e_{i-1}$. Then, by Theorem \ref{tcnn},
$$\Hom_{\Lambda(n,\overline{m})}(\Omega(M_l),M_l)=\langle f_t:t\in \{l+1,l+2,\cdots,i-1\}\rangle,$$
where $f_t$ is the canonical homomorphism $\Omega(M_l) \twoheadrightarrow S(t) \xhookrightarrow{} M_l$. We will prove that $f_t\in \mathcal{P}_{\Lambda(n,\overline{m})}(\Omega(M_l),M_l)$ for all $t\in \{l+1,l+2,\cdots,i-1\}$. First, observe that
\begin{equation}\label{efl1}
f_{l+1}=\psi_{l+1} \varphi_{l+1},
\end{equation}
where $\varphi_{l+1}$ and $\psi_{l+1}$ are the canonical homomorphisms
$$\varphi_{l+1}: \Omega(M_l) \twoheadrightarrow M[x_{l+1,l}A_l^{m_0-1}] \xhookrightarrow{} P(l),$$
$$\psi_{l+1}: P(l) \twoheadrightarrow M[\delta_t^{-(m_{l+1}-1)}\alpha_1] \xhookrightarrow{} M_l.$$
Now, suppose that $t\in \{l+2,\cdots,i-1\}$. Def\mbox ine $\varphi_t$ and $\psi_t$ to be the canonical homomorphisms
$$\varphi_t: \Omega(M_l) \twoheadrightarrow M[x_{t,t-1}A_{t-1}^{m_0-1}\delta_{t-1}^{-1}] \xhookrightarrow{} P(t-1),$$
$$\psi_t: P(t-1) \twoheadrightarrow M[\delta_{t-1}^{-(m_t-1)}\alpha_{t-1}] \xhookrightarrow{} M_l.$$
From the forms of the homomorphisms $f_t$, $\varphi_t$ and $\psi_t$, we have
\begin{equation}\label{eft}
f_t=\psi_t\varphi_t-f_{t-1}.
\end{equation}
Since $\psi_t\varphi_t\in \mathcal{P}_{\Lambda(n,\overline{m})}(\Omega(M_l),M_l)$ for all $t\in \{l+1,\cdots,i-1\}$, it follows from \eqref{efl1} and \eqref{eft} that $f_t\in \mathcal{P}_{\Lambda(n,\overline{m})}(\Omega(M_l),M_l)$ for all $t\in \{l+1,\cdots,i-1\}$. In consequence,
$$\Ext_{\Lambda(n,\overline{m})}^1(M_l,M_l)\cong \shom_{\Lambda(n,\overline{m})}(\Omega(M_l),M_l)=0.$$

\item[\textbf{Case 3}] $l=0$.\\
As in Case 2, a similar calculation yields $M_0=M[y_0x_{i,n}]$, $\Omega(M_0)=M[S]$ with
$$
S=x_{0,i-1}A_{i-1}^{m_0-1} \cdot x_{i-1,i-2}A_{i-2}^{m_0-1} \cdot\ldots \cdot x_{1,0}A_0^{m_0-1}, 
$$
and $\Hom_{\Lambda(n,\overline{m})}(\Omega(M_l),M_l)$ is generated by the canonical homomorphisms $f_t: \Omega(M_0) \twoheadrightarrow S(t) \xhookrightarrow{} M_0$ for $t\in \{0,1\ldots,i-1\}$, and $f_t\in \mathcal{P}_{\Lambda(n,\overline{m})}(\Omega(M_0),M_0)$ for $t\in \{1,\ldots,i-1\}$. To show that $f_0\notin \mathcal{P}_{\Lambda(n,\overline{m})}(\Omega(M_0),M_0)$, suppose there exist a projective module $P$ and homomorphisms $\varphi: \Omega(M_0) \rightarrow P$ and $\psi: P \rightarrow M_0$ such that $\psi \varphi$ maps the quotient $S(0)$ of $\Omega(M_0)$ to the corresponding submodule of $M_0$ of the same form. By the Krull-Schmidt-Azumaya Theorem, we can assume $P$ is indecomposable. The structure of the radical series of $M_0$, implies $P\cong P(0)$. However, due to the radical series of $P(0)$, the image of $S(0)$ under $\varphi$ and the preimage of $S(0)$ under $\psi$ lie in dif\mbox ferent radical series, which is a contradiction. Thus, $S(0)$ is in the kernel of every homomorphism in $\mathcal{P}_{\Lambda(n,\overline{m})}(\Omega(M_0),M_0)$, implying that $f_0\notin \mathcal{P}_{\Lambda(n,\overline{m})}(\Omega(M_0),M_0)$. In summary, $f_t$ is in $\mathcal{P}_{\Lambda(n,\overline{m})}(\Omega(M_0),M_0)$ for all $t\neq0$. Therefore,
$$\dim_{\Bbbk}\Ext_{\Lambda(n,\overline{m})}^1(M_0,M_0)=\dim_{\Bbbk}\shom_{\Lambda(n,\overline{m})}(\Omega(M_0),M_0)=1.$$
\end{description}
\end{proof}

\begin{theo}\label{tudrgbt} Suppose that at least one vertex of $\mathcal{W}_{n,\overline{m}}$ has multiplicity $1$. If $M$ is a periodic string $\Lambda(n,\overline{m})$-module such that $\send_{\Lambda(n,\overline{m})}(M)\cong \Bbbk$, then
\begin{equation*}
R(\Lambda(n,\overline{m}),M)= \begin{cases}
\Bbbk & \text{if $d_M<n$,}\\
\Bbbk \llbracket x \rrbracket & \text{if $d_M=n$.}
\end{cases}
\end{equation*}
\end{theo}

\begin{proof}
If $d_M<|E|-1$, by Theorem \ref{text} we have that $\Ext_{\Lambda(n,\overline{m})}^1(M,M)=0$ and hence $R(\Lambda,M)=\Bbbk$.

For the case $d_M=n$, we will use Theorem \ref{anote}. By Theorem \ref{text}, $\dim_{\Bbbk} \Ext_{\Lambda(n,\overline{m})}^1(M,M)=1$, implying that $R(\Lambda(n,\overline{m}),M)$ is a quotient of $\Bbbk \llbracket x \rrbracket$ (see \cite[Remark 9]{bw}). Given Theorem \ref{trsyz} and Lemma \ref{lsmmp}, its suf\mbox f\mbox ices to consider $M[y_0x_{i,n}]$. Def\mbox ine $\mathcal{L}:=\{V_0,V_1,\ldots\}$, where $V_0=M_0$ and $V_l=N_{nl}$ for $l>0$ (see Remark \ref{rd} for the def\mbox inition of the string modules $N_p$). Then, $V_0=M[y_0x_{i,n}]$ and $V_l=M[S_l]$ for $l>0$, with $S_l$ def\mbox ined recursively by $S_l=y_0x_{i,n}\alpha_nS_{l-1}$ and $S_0:=y_0x_{i,n}$. For $l>0$, let $\iota_l$ be the natural inclusion $V_{l-1}\xhookrightarrow{} V_l$ induced by the factorization $S_l=y_0x_{i,n}\alpha_nS_{l-1}$ and $\epsilon_l$ the canonical projection $V_l\twoheadrightarrow V_{l-1}$ induced by $S_l=S_{l-1}\alpha_ny_0x_{i,n}$. Set $\sigma_l:=\iota_l \epsilon_l$. To verify the hypothesis of Theorem \ref{anote}, f\mbox irst observe that $\ker \sigma_l=V_0$. By induction on $t\in \{1,\ldots, l\}$ and the form of the strings $S_l$, we have $\im(\sigma_l^t)=V_{l-t}$, implying $\im(\sigma_l^l)=V_0$. Since $\mathcal{L}$ is inf\mbox inite, Theorem \ref{anote} yields $R(\Lambda(n,\overline{m}),M)=\Bbbk \llbracket x \rrbracket$.
\end{proof}

\subsubsection{Second case of the classif\mbox ication}\label{subsec:two}

For this case, we assume that all vertices of $\mathcal{W}_{n,\overline{m}}$ have multiplicity greater than 1. While some proofs may follow similar reasoning as in the previous case, we will streamline the arguments to avoid redundancy.

All the projective indecomposable modules $P(j)$, for $j\in\{0,\ldots,n\}$, have the following radical series
\vspace{-8mm}

\[
P(j):\ \tabbedCenterstack{& \mbox{$S(j)$} &&& \\
\mbox{$S(j)$} \ \ && \ \ \mbox{$S(j+1)$} && \\
&&&& \\
\vdots && \ \ \vdots && \\
&& \ \ \mbox{$S(j)$} && \\
\mbox{$S(j)$} \ \ &&&& \\
\ \ && \mbox{$S(j+1)$} && \\
&&&& \\
\vdots \ \ && \ \ \vdots && \\
&&&& \\
\mbox{$S(j)$} \ \ && \ \ \mbox{$S(j-1)$} && \\
& \mbox{$S(j)$} &&& }
\]

\noindent
As there are no simple uniserial $\Lambda(n,\overline{m})$-modules, modules on the boundaries of exceptional tubes of $\leftindex_s \Gamma_{\Lambda(n,\overline{m})}$ are maximal uniserial submodules of projective indecomposable modules. These modules are the string modules $M[x_{j+1,j}A_j^{m_0-1}]$ and $M[\delta_j^{m_{j+1}-1}]$, where $j\in \{0,\ldots,n\}$. To describe modules $M$ such that $d_M\leq n$, we def\mbox ine the string
\begin{equation*}
z_j:= \begin{cases}
\delta_n^{-(m_{n+1}-1)} & \text{if $j=n$,}\\
\left(\delta_j^{-(m_{j+1}-1)}\alpha_j\right) \cdots \left(\delta_{n-1}^{-(m_n-1)}\alpha_{n-1}\right) \delta_n^{-(m_{n+1}-1)} & \text{if $j\in \{0,\ldots,n-1\}$.}
\end{cases}
\end{equation*}

\begin{lemm}\label{lstr22} Let $M$ be a periodic string $\Lambda(n,\overline{m})$-module. Then, $d_M\leq n$ if and only if $M$ is in the $\Omega$-orbit of a string module $M[z_j]$, for some $j\in \{0,\ldots,n\}$.
\end{lemm}

\begin{proof}
It is enough to take the f\mbox irst $n+1$ vertices of a diagonal $\mathfrak{C}$ of some exceptional tube of $\leftindex_s\Gamma_{\Lambda(n,\overline{m})}$. Since $M[z_n]$ is on the boundary of an exceptional tube, there exists a sectional path
\[
\xymatrix{
M_n \ar[r] & M_{n-1} \ar[r] & \cdots \ar[r] & M_0
}
\]
in $\leftindex_s\Gamma_{\Lambda(n,\overline{m})}$, with $M_n:=M[z_n]$. For each $j\in \{1,\ldots,n-1\}$, $M_j=(M_{j-1})_h$, implying $M_j=M[z_j]$.
\end{proof}

\begin{remark}\label{rst2np} \rm The whole diagonal of $\mathfrak{C}$ starting at $M_n$ is the following
\[
\xymatrix{
M_n \ar[r] & \cdots \ar[r] & M_0 \ar[r] & N_0 \ar[r] & N_1 \ar[r] & \cdots
}
\]
where the modules $M_j$, for $0\leq j\leq n$, are as in the proof of Lemma \ref{lstr22}, $N_0:=(M_0)_h$ and $N_p:=(N_{p-1})_h$ for all $p>0$. Let $C_p$ be the string such that $N_0=M[C_p]$ and let $q$ and $r$ be the unique integers such that $p=(n+1)q+r$ and $0\leq r\leq n$. Thus, we have $C_p=z_{n-r}(\alpha_nz_0)^{q+1}$.
\end{remark}

\begin{theo}\label{tstables2} Suppose that all the vertices in $\mathcal{W}_{n,\overline{m}}$ have multiplicity greater than $1$. Let $M$ be a string periodic module. The stable endomorphism ring of $M$ is isomorphic to $\Bbbk$ if and only if $d_M\leq n$.
\end{theo}

\begin{proof}
As $\dim_{\Bbbk} \send_{\Lambda(n,\overline{m})}(M)=\dim_{\Bbbk} \send_{\Lambda(n,\overline{m})}(\Omega(M))$ (by Lemma \ref{lstr22} and Remark \ref{rst2np}) it suf\mbox f\mbox ices to show that $\send_{\Lambda(n,\overline{m})}(M_j)\cong \Bbbk$ for $j\in \{0,\ldots,n\}$ and $\send_{\Lambda(n,\overline{m})}(N_p)\not\cong \Bbbk$ for $p\in \{0,1,\ldots\}$.

We f\mbox irst show that $M_j$ has stable endomorphism ring $\Bbbk$ for $j\in \{0,\ldots,n\}$. A canonical endomorphism $f_t:M_n \xhookrightarrow{} M[C] \twoheadrightarrow M_n$ distinct from the identity exists if and only if $C=\delta_t^l$ for some $t\in \{j,\ldots,n\}$ and $l\in \{1,\ldots,m_{l+1}-2\}$. Here, $C:=\delta_t^0:=e_t$. We consider the following cases based on the value of $t$.
\begin{description}
\item[\textbf{Case 1}] For $t=j$, we have that $f_j=\psi_j\varphi_j$, where $\varphi_j$ and $\psi_j$ are the canonical homomorphisms $M_j\twoheadrightarrow M[\delta_j^{-(m_{j+1}-1)}] \xhookrightarrow{} P(j)$ and $P(j)\twoheadrightarrow M[\delta_j^{-(l+1)}] \xhookrightarrow{} M_j$, respectively.
\item[\textbf{Case 2}] For $t\neq j$, we have that $f_t=\psi_t\varphi_t$, where $\varphi_t$ and $\psi_t$ are the canonical homomorphisms $M_j\twoheadrightarrow M[\alpha_{t-1}\delta_t^{-(m_{j+1}-1)}] \xhookrightarrow{} P(t)$ and $P(t)\twoheadrightarrow M[\delta_j^{-(l+1)}] \xhookrightarrow{} M_j$, respectively.
\end{description}
In both cases, $f_t$ factors through a projective module, implying $\send_{\Lambda(n,\overline{m})}(M_j)\cong \Bbbk$ for every $j\in \{0,\ldots,n\}$. On the other hand, by Remark \ref{rst2np}, $C_p=z_{n-r}\cdot (\alpha_n z_0)^{q+1}$. Additionally, the factorization
$$C_p=\left(z_{n-r}(\alpha_n z_0)^q(\delta_0^{-(m_1-1)}\alpha_0)\cdots (\delta_{n-p-1}^{-(m_{n-p}-1)}\alpha_{n-p-1})\right)\cdot z_{n-r},$$
induces a canonical endomorphism $N_p\twoheadrightarrow M[z_{n-r}] \xhookrightarrow{} N_p$ that does not factor through a projective module. Thus, $\send_{\Lambda(n,\overline{m})}(N_p)\not\cong \Bbbk$ for every $p\in \{0,1,\ldots\}$.
\end{proof}

\begin{theo}\label{text22} Suppose that all the vertices in $\mathcal{W}_{n,\overline{m}}$ have multiplicity greater than $1$. Let $M$ be a string periodic $\Lambda(n,\overline{m})$-module with stable endomorphism ring isomorphic to $\Bbbk$. Then,
\begin{equation*}
\dim_{\Bbbk}\Ext_{\Lambda(n,\overline{m})}^1(M,M)= \begin{cases}
0 & \text{if $d_M<|E|-1$,}\\
1 & \text{if $d_M=|E|-1$.}
\end{cases}
\end{equation*}
\end{theo}

\begin{proof}
As in the proof of Theorem \ref{tstables2}, it is enough to consider the modules $M_j$. Note that
\begin{equation*}
\Omega(M_j)= \begin{cases}
M[x_{0,n}A_n^{m_0-1}] & \text{if $j=n$,}\\
M[x_{0,n}A_n^{m_0-1}](\delta_n^{-1}x_{n,n-1}A_{n-1}^{m_0-1})\cdots(\delta_{j+1}^{-1}x_{j+1,j}A_j^{m_0-1}) & \text{if $j\neq n$.}
\end{cases}
\end{equation*}
For $j=n$, there are no canonical homomorphisms from $\Omega(M_n)$ to $M_n$ and then
$$
\dim_{\Bbbk}\Ext_{\Lambda(n,\overline{m})}^1(M_n,M_n)=\dim_{\Bbbk} \shom_{\Lambda(n,\overline{m})}(\Omega(M_n),M_n)=0.
$$
For $j\neq n$, $\Hom_{\Lambda(n,\overline{m})}(\Omega(M_n),M_n)$ is spanned by canonical homomorphisms $f_t:\Omega(M_j)\twoheadrightarrow S(t)\xhookrightarrow{} M_j$, where $t\in \{j+1,\ldots,n\}$ if $j\neq 0$ and $t\in \{0,1,\ldots,n\}$ if $j=0$. For $t=j+1$, $f_{j+1}=\psi_{j+1}\varphi_{j+1}$, with $\varphi_{j+1}: \Omega(M_j)\twoheadrightarrow M[x_{j+1,j}A_j^{m_0-1}] \xhookrightarrow{} P(j)$ and $\psi_{j+1}: P(j)\twoheadrightarrow M[\delta_j^{-(m_{j+1}-1)}\alpha_j] \xhookrightarrow{} M_j$ canonical homomorphisms. Thus, $f_{j+1}\in \mathcal{P}_{\Lambda(n,\overline{m})}(\Omega(M_j),M_j)$. For $t\neq j+1$ and $t\neq0$, $f_t=\psi_t\varphi_t-f_{t-1}$, where $\varphi_t: \Omega(M_j)\twoheadrightarrow M[x_{t,t-1}A_{t-1}^{m_0-1}\delta_{t-1}^{-1}] \xhookrightarrow{} P(t-1)$ and $\psi_t:P(t-1)\twoheadrightarrow M[\delta_{t-1}^{-(m_t-1)}\alpha_{t-1}] \xhookrightarrow{} M_j$ canonical homomorphisms. Therefore, $f_t\in \mathcal{P}_{\Lambda(n,\overline{m})}(\Omega(M_j),M_j)$. Consequently, $\Ext_{\Lambda(n,\overline{m})}^1(M_j,M_j)=0$ for $j\neq 0$.

For $j=0$, there are no canonical homomorphisms $\Omega(M_0)\twoheadrightarrow M[C]\xhookrightarrow{} M_0$ with $C=e_0D$, for some string $D$. Therefore,  $f_0\notin \mathcal{P}_{\Lambda(n,\overline{m})}(\Omega(M_0),M_0)$, implying $$
\dim_{\Bbbk}\Ext_{\Lambda(n,\overline{m})}^1(M_0,M_0)=\dim_{\Bbbk} \shom_{\Lambda(n,\overline{m})}(\Omega(M_0),M_0)=1.
$$
\end{proof}

\begin{theo}\label{tudr22} Suppose that all the vertices in $\mathcal{W}_{n,\overline{m}}$ have multiplicity greater than $1$. Let $M$ be a string periodic $\Lambda(n,\overline{m})$-module with stable endomorphism ring isomorphic to $\Bbbk$. Then,
\begin{equation*}
R(\Lambda(n,\overline{m}),M)= \begin{cases}
\Bbbk & \text{if $d_M<n$,}\\
\Bbbk \llbracket x \rrbracket & \text{if $d_M=n$.}
\end{cases}
\end{equation*}
\end{theo}

\begin{proof}
If $d_M<n$, then $\Ext_{\Lambda(n,\overline{m})}^1(M,M)=0$, {\it a fortiori}, $R(\Lambda(n,\overline{m}),M)=\Bbbk$, by Theorem \ref{text22}. Suppose that $d_M=n$. Since $\Ext_{\Lambda(n,\overline{m})}^1(M,M)\cong \Bbbk$, the universal deformation ring $R(\Lambda(n,\overline{m}),M)$ is a quotient of $\Bbbk \llbracket x \rrbracket$. We def\mbox ine the sequence $\mathcal{L}=\{M[T_0],M[T_1],\ldots\}$, where $T_0=z_0$ and $T_l=z_0\alpha_nT_{l-1}$ for $l>0$. By the def\mbox inition of $T_l$, we have a natural inclusion $\iota_l:M[T_{l-1}]\xhookrightarrow{} M[T_l]$, and the alternative factorization $T_l=T_{l-1}\alpha_nz_0$ induces the canonical epimorphism $\epsilon:M[S_l]\twoheadrightarrow M[T_{l-1}]$. Set $\sigma_l:=\iota_l\epsilon_l$. We have $\ker \sigma_l=T_0$, and by the form of the strings $T_l$, it follows that $\im(\sigma_l^t)=M[T_{l-t}]$ for $t\in \{1,\ldots, l\}$. In particular, $\im(\sigma_l^l)=M[T_0]$. Since $\mathcal{L}$ is inf\mbox inite, by Theorem \ref{anote} we conclude that $R(\Lambda(n,\overline{m}),M)=\Bbbk \llbracket x \rrbracket$.
\end{proof}

We summarize the classif\mbox ication results from subsections \ref{subsec:one} and \ref{subsec:two} in the following theorem.

\begin{theo}\label{udr} Let $\mathcal{G}=(V,E,\mathfrak{m},\mathfrak{o})$ be a generalized Brauer tree and $\Lambda_{\mathcal{G}}$ its associated Brauer graph algebra. Suppose that $M$ is a periodic string $\Lambda_{\mathcal{G}}$-module. Then, $M$ has stable endomorphism ring isomorphic to $\Bbbk$ if and only if $d_M\leq |E|-1$. Moreover, if $M$ has stable endomorphism ring isomorphic to $\Bbbk$, then
\begin{equation*}
R(\Lambda_{\mathcal{G}},M)= \begin{cases}
\Bbbk & \text{if $d_M<|E|-1$,}\\
\Bbbk \llbracket x \rrbracket & \text{if $d_M=|E|-1$.}
\end{cases}
\end{equation*}
\end{theo}

\begin{proof}
By Lemma \ref{lbdes}, $\Lambda_{\mathcal{G}}$ is derived equivalent to $\Lambda(n,\overline{m})$, where $n=|E|-1$ and $\overline{m}=(m_0,\ldots,m_i)$ is the vector of multiplicities greater than $1$ of $\mathcal{G}$. Then from Theorem \ref{tderst}, $\Lambda_{\mathcal{G}}$ and $\Lambda(n,\overline{m})$ are stably equivalent of Morita type. Moreover, Lemma \ref{lcgbt} implies $d_M=d_{\Omega(M)}$. Since universal deformation rings are invariant under stable equivalence of Morita type (see Theorem \ref{tedstm}), the statement holds for $\Lambda_{\mathcal{G}}$ if and only if it holds for $\Lambda(n,\overline{m})$. Thus, the statement follows from Theorems \ref{tmek}, \ref{tudrgbt}, \ref{tstables2} and \ref{tudr22}.
\end{proof}

\subsection{An application: over Standard-Koszul algebras}

In this subsection, we apply the results from this section to compute universal deformation rings of modules in exceptional tubes over standard Koszul symmetric special biserial algebras. First, we present some preliminaries.

As shown in \cite{mgy}, A. Magyar studies the connection between the standard Koszul and quasi-Koszul properties. In particular, Magyar characterizes certain symmetric special biserial algebras that are standard Koszul. Since symmetric special biserial algebras coincide with Brauer graph algebras, and the algebras in Magyar's characterization are generalized Brauer tree algebras, we can apply Theorem  \ref{udr} to compute their universal deformation rings.

\begin{defi}[\cite{mgy}]\label{dmoddel} Let $\Lambda$ be a f\mbox inite-dimensional algebra and $\{e_1\ldots,e_n\}$ a complete family of primitive orthogonal idempotents of $\Lambda$. Def\mbox ine $\varepsilon_i:=e_i+\cdots+e_n$ for $i\in \{1,\ldots,n\}$ and set $\varepsilon_0:=0$. The \emph{standard right module} $\Delta(i)$ is def\mbox ined to be $e_i\Lambda/e_i\Lambda \varepsilon_{i+1} \Lambda$ and the \emph{standard left module} $\overline{\Delta}^{\circ}(i)$ is def\mbox ined to be $\Lambda e_i/\Lambda \varepsilon_{i+1} \rad(\Lambda) e_i$
\end{defi}

\begin{defi}[\cite{mgy}]\label{dclassci} Let $\Lambda$ be a f\mbox inite-dimensional algebra. A module $M$ of$\mod\hspace{-1.7mm} -\hspace{-.2mm} \Lambda$ is in the class $\mathcal{C}_{\Lambda}$ if $\Omega_j\subseteq \rad(P_{j-1})$ holds for every $j\geq 0$, where
\[
\xymatrix{
\cdots \ar[r]^{\alpha_2} & P_1 \ar[r]^{\alpha_1} & P_0 \ar[r]^{\alpha_0} & M
}
\]
is a minimal projective resolution of $M$ and $\Omega_j:=\ker(\alpha_{j-1})$ for $j\in \{1,2,\ldots\}$.
\end{defi}

\begin{defi}[\cite{mgy}]\label{dkszl} A f\mbox inite-dimensional algebra $\Lambda$ is said to be a \emph{standard Koszul algebra} if $\Delta(i)$ is in $\mathcal{C}_{\Lambda}$ and $\overline{\Delta}^{\circ}(i)$ is in $\mathcal{C}_{\Lambda^{op}}$ for all $i\in \{1,\ldots,n\}$.
\end{defi}

In the following theorem, we present Magyar's characterization for standard Koszul symmetric special biserial algebras.

\begin{theo}[\cite{mgy}]\label{tmgyr} Let $\Lambda=\Bbbk Q/I$ be a symmetric special biserial algebra. Then, $\Lambda$ is a standard Koszul algebra if and only if $Q$ and $I$ have the form
$$
Q:\xymatrix{
0 \ \ar@(ul,dl)_{\delta} \ar@<-0.5ex>[rr]_{\alpha_1} & & \ 1 \ \ar@<-0.5ex>[ll]_{\beta_1} \ \ar@<-0.5ex>[rr]_{\alpha_2} & & \ \ar@<-0.5ex>[ll]_{\beta_2} \cdots \ \ar@<-0.5ex>[rr]_{\alpha_{n-1}} & & \ n-1 \ \ar@<-0.5ex>[ll]_{\beta_{n-1}} \ \ar@<-0.5ex>[rr]_{\alpha_n} & & \ n \ \ar@<-0.5ex>[ll]_{\beta_n} \ar@(dr,ur)_{\gamma}}
$$
$$
I=\langle \delta^l-\alpha_1\beta_1,\lambda \gamma^m-\beta_n\alpha_n,\beta_i\alpha_i-\alpha_{i+1}\beta_i,\delta\alpha_1, \beta_1\delta,\alpha_n\gamma,\gamma\beta_n,\alpha_i\alpha_{i+1},\beta_{i+1}\beta_i:1\leq i\leq n-1\rangle
$$
for some $n,l,m\in \mathbb{Z}^+$ and $\lambda \in \Bbbk\backslash \{0\}$, with $l\geq 2$ and $m\geq 2$.
\end{theo}

From now on, for each $\lambda \in \Bbbk\backslash\{0\}$, we write $\Lambda_{\lambda}$ to denote the bound quiver algebra presented in Theorem \ref{tmgyr}. Observe that the radical series of the projective indecomposable modules over $\Lambda_{\lambda}$ are the following.

$$P(0):\ \tabbedCenterstack{
& \mbox{$S(0)$} & \\
&& \\
\mbox{$S(0)$} && \\
&& \\
\vdots \ && \ \ \mbox{$S(1)$} \\
&& \\
\mbox{$S(0)$} && \ \ \\
&&\\
& \mbox{$S(0)$} & }\ , \ \ \ \
P(n):\ \tabbedCenterstack{
& \mbox{$S(n)$} & \\
&& \\
\mbox{$S(n)$} && \\
&& \\
\vdots \ && \mbox{$S(n-1)$} \\
&& \\
\mbox{$S(n)$} && \ \ \\
&&\\
& \mbox{$S(n)$} &
}\ ,$$
\vspace{3mm}

$$\mbox{ and }\ \
P(j):\ \tabbedCenterstack{
& \mbox{$S(j)$} & \\
&& \\
\mbox{$S(j-1)$} && \mbox{$S(j+1)$}\\
&& \\
& \mbox{$S(j)$} &
}\ , \mbox{ for } j\in \{1,\ldots,n-1\}.
$$
\vspace{3mm}

Note that $\Lambda_{\lambda}$ it is also a Brauer graph algebra. The following theorem shows that $\Lambda_{\lambda}$ is derived equivalent to a star algebra $\mathcal{W}_{n,\overline{m}}$ as def\mbox ined in Subsection \ref{sedis}.

\begin{figure}[h!]
\begin{center}
\begin{tikzpicture}[square/.style={regular polygon,regular polygon sides=4,minimum width=14mm,inner sep=0pt}]
\tikzset{vnode/.style={draw,thick,circle,minimum width=9mm,inner sep=0pt}};
    \node at (0,0) [square,draw] (zeta0) {$\zeta_0$};
    \node at (2,0) [vnode,draw] (zeta1) {$\zeta_1$};
    \node at (4,0) (cdots) {\textbf{$\cdots$}};
    \node at (6,0) [vnode,draw] (zetan) {$\zeta_n$};
    \node at (8,0) [square,draw] (zetan1) {$\zeta_{n+1}$};
    \draw (zeta0) -- (zeta1) -- (cdots) -- (zetan) -- (zetan1);
\end{tikzpicture}
\end{center}
\caption{The star $\mathcal{W}_{n,\overline{m}}$ associated to $\Lambda_{\lambda}$}\label{fig3}
\end{figure}
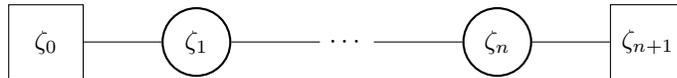

\begin{lemm}\label{lsteq} The algebra $\Lambda_{\lambda}$ is derived equivalent to a generalized Brauer tree algebra whose underlying Brauer graph is a star with $n+1$ edges (see Figure \ref{fig3}). This star has two vertices, $\zeta_0$ and $\zeta_{n+1}$, with multiplicities greater than 1 (namely, $l$ and $m$, respectively) and all other vertices have multiplicity $1$. 
\end{lemm}

\begin{proof} By \cite[Lemma 4.2]{kean1}, it follows that $\Lambda_{\lambda}$ is isomorphic to $\Lambda_1=\Bbbk Q/I_1$, where $I_1=\langle \delta^l-\alpha_1\beta_1, \gamma^m-\beta_n\alpha_n,\beta_i\alpha_i-\alpha_{i+1}\beta_i,\delta\alpha_1, \beta_1\delta,\alpha_n\gamma,\gamma\beta_n,\alpha_i\alpha_{i+1},\beta_{i+1}\beta_i:1\leq i\leq n-1\rangle$. Then, $\Lambda_1$ is the Brauer graph algebra whose underlying Brauer graph is the generalized Brauer tree in Figure \ref{fig3}.

By Theorem \ref{lbdes}, $\Lambda_1$ is derived equivalent to a Brauer graph algebra whose underlying Brauer graph is a star with $n+1$ edges and two vertices of multiplicities greater than $1$.
\end{proof}

\begin{coro}\label{cetksz} There are two exceptional tubes in the stable Auslander-Reiten quiver of $\Lambda_{\lambda}$ and they are of rank $2n+2$. Moreover, all the modules on their boundaries coincide with the maximal uniserial submodules of the projective indecomposable $\Lambda_{\lambda}$-modules.
\end{coro}

\begin{proof}
By Lemmas \ref{lsteq} and \ref{lbdes}, there is a rank-preserving bijection between the exceptional tubes in the stable Auslander-Reiten quivers of $\Lambda_{\lambda}$ and $\Lambda(n,\overline{m})$. As there are no simple uniserial $\Lambda_{\lambda}$-modules, the modules on the boundaries of the two exceptional tubes in $\leftindex_s\Gamma_{\Lambda_{\lambda}}$ are maximal uniserial submodules of projective indecomposable $\Lambda_{\lambda}$-modules.
\end{proof}

For a more explicit description, the modules on the boundaries of the exceptional tubes of $\leftindex_s\Gamma_{\Lambda_{\lambda}}$ are the string modules $M[C]$, where $C\in \{\delta^{l-1},\gamma^{m-1},\alpha_j,\beta_j:j=1,\ldots,n\}$.

Combining Theorem \ref{udr}, its proof, and Lemma \ref{lsteq}, we obtain the following result.

\begin{theo}\label{tksz} Let $\Lambda_{\lambda}=\Bbbk Q/I$ be a standard Koszul symmetric special biserial algebra and $M$ be a periodic string module over $\Lambda_{\lambda}$. Then, $M$ has stable endomorphism ring isomorphic to $\Bbbk$ if and only if $d_M\leq |Q_0|-1$. Moreover, if $M$ has stable endomorphism ring isomorphic to $\Bbbk$, then
\begin{equation*}
R(\Lambda_{\lambda},M)= \begin{cases}
\Bbbk & \text{if $d_M<|Q_0|-1$,}\\
\Bbbk \llbracket x \rrbracket & \text{if $d_M=|Q_0|-1$.}
\end{cases}
\end{equation*}
\end{theo}

\begin{remark} 
\rm To explicitly describe the modules in Theorem \ref{tksz}, consider the sectional path in $\leftindex_s\Gamma_{\Lambda_{\lambda}}$ starting at $M[\rho^{m-1}]$. A periodic string $\Lambda_{\lambda}$-module $M$ has stable endomorphism ring isomorphic to $\Bbbk$ if and only if it belongs to the $\Omega$-orbit of a string module corresponding to a string in one of the following cases.
\begin{itemize}
   \item[$i)$] $\gamma^{m-1}$, $\gamma^{m-1}\alpha_n^{-1}\beta_{n-1}$, ..., $\gamma^{m-1}\alpha_n^{-1}\beta_{n-1}\cdots \alpha_2^{-1}\beta_1$,

$\gamma^{m-1}\alpha_n^{-1}\beta_{n-1}\cdots \alpha_2^{-1}\beta_1\delta^{-1}\alpha_1$,

$\gamma^{m-1}\alpha_n^{-1}\beta_{n-1}\cdots \alpha_2^{-1}\beta_1\delta^{-1}\alpha_1\beta_2^{-1}\alpha_3$, ...,

$\gamma^{m-1}\alpha_n^{-1}\beta_{n-1}\cdots \alpha_2^{-1}\beta_1\delta^{-1}\alpha_1\beta_2^{-1}\alpha_3\cdots \beta_{n-2}^{-1}\alpha_{n-1}$,

for $n$ even.
   \item[$ii)$] $\gamma^{m-1}$, $\gamma^{m-1}\alpha_n^{-1}\beta_{n-1}$, ..., $\gamma^{m-1}\alpha_n^{-1}\beta_{n-1}\cdots \alpha_3^{-1}\beta_2$,

$\gamma^{m-1}\alpha_n^{-1}\beta_{n-1}\cdots \alpha_3^{-1}\beta_2\alpha_1^{-1}\delta^{l-1}$,

$\gamma^{m-1}\alpha_n^{-1}\beta_{n-1}\cdots \alpha_3^{-1}\beta_2\alpha_1^{-1}\delta^{l-1}\beta_1^{-1}\alpha_2$, ...,

$\gamma^{m-1}\alpha_n^{-1}\beta_{n-1}\cdots \alpha_3^{-1}\beta_2\alpha_1^{-1}\delta^{l-1}\beta_1^{-1}\alpha_2 \cdots \beta_{n-2}^{-1}\alpha_{n-1}$,

for $n$ odd.

\end{itemize}
\end{remark}

\section{Non-periodic string modules and its Universal Deformation Rings: a special case of classif\mbox ication}\label{nonperiodiccase}

Let $\mathcal{W}_{n,\overline{m}}$ be a star as def\mbox ined in Subsection \ref{sedis}. For a non-periodic simple module $S(t)$, let $\mathfrak{C}_{S(t)}$ denote the component of $\leftindex_s\Gamma_{\Lambda}$ containing $S(t)$. In this section we give necessary and suf\mbox f\mbox icient conditions for $\mathfrak{C}_{S(t)}$ to be $\Omega$-stable. Moreover, assuming that $\mathfrak{C}_{S(t)}$ is $\Omega$-stable, we compute and classify the universal deformation rings of modules $M$ in $\mathfrak{C}_{S(t)}$, provided that $\send_{\Lambda(n,\overline{m})}(M)\cong \Bbbk$.

\subsection{$\Omega$-stable components containing simple modules}

Let $\mathcal{W}_{n,\overline{m}}$ be a star, where $\overline{m}=(m_0,\ldots,m_i)$. Lemma \ref{le} provides necessary and suf\mbox f\mbox icient conditions for a component containing a non-periodic simple module to be  $\Omega$-stable.

\begin{lemm}\label{le} Let $S(t)$ be a non-periodic simple module in\hspace{-2mm} $\mod\hspace{-1.5mm} -\hspace{-.3mm} \Lambda(n,\overline{m})$. Then, $\mathfrak{C}_{S(t)}$ is $\Omega$-stable if and only if at least one of the following conditions holds.
\begin{enumerate}
    \item[$i)$] $m_{t+1}=2$;
    \item[$ii)$] $i=1$, $t=0$, and $m_0=2$.
\end{enumerate}
\end{lemm}

\begin{proof} Note that $\Omega(S(t))=\rad(P(i))$. Thus, $\mathfrak{C}_{S(t)}$ is $\Omega$-stable if and only if $\rad(P(t))$ belongs to $\mathfrak{C}_{S(t)}$. Now, suppose that $\mathfrak{C}_{S(t)}$ is $\Omega$-stable. Since $t\in \{0,\ldots,i-1\}$, condition $ii)$ is equivalent to $i=1$ and $m_0=2$. We claim that $m_{t+1}\neq 2$ implies $i=1$ and $m_0=2$. Indeed, by Theorem \ref{tstcsrad}, there is a path
{\footnotesize
\begin{center}
$p:$\ \
\begin{tikzpicture}[square/.style={regular polygon,regular polygon sides=4,minimum width=14mm,inner sep=0pt}]
\tikzset{vnode/.style={draw,thick,circle,minimum width=9mm,inner sep=0pt}};
    \node at (0,0) (u0) {$v_0$};
    \node at (2,0) (v1) {$v_1$};
    \node at (4,0) (cdots) {\textbf{$\cdots$}};
    \node at (6,0) (vt1) {$v_{h-1}$};
    \node at (8,0) (vt) {$v_h$\ .};
    \draw (u0) edge node [above] {$a_1=t$} (v1);
    \draw (v1) edge node [above] {$a_2$} (cdots);
    \draw (cdots) edge node [above] {$a_{h-1}$} (vt1);
    \draw (vt1) edge node [above] {$a_h=t$} (vt);
\end{tikzpicture}
\end{center}
}
\hspace{-6.8mm} The path $p$ must have even length and satisfy the following properties, for all $k\in \{1,\ldots,h-1\}$:
\begin{description}
\item[\textbf{Property 1}] $m(v_k)=1$ if $a_k\neq a_{k+1}$  and $m(v_k)=2$ if $a_k=a_{k+1}$.
\item[\textbf{Property 2}] $a_k$ and $a_{k+1}$ are the unique non-exceptional edges incident to $v_k$.
\end{description}
Since $m_{t+1}>2$ and $m_0>1$, the f\mbox irst condition forces $a_1=\cdots=a_h=t$ and $v_k\neq \zeta_{t+1}$ for every $k$. Thus, $p$ reduces to the path
{\footnotesize
\begin{center}
\begin{tikzpicture}[square/.style={regular polygon,regular polygon sides=4,minimum width=14mm,inner sep=0pt}]
\tikzset{vnode/.style={draw,thick,circle,minimum width=9mm,inner sep=0pt}};
    \node at (0,0) (u0) {$\zeta_{t+1}$};
    \node at (2,0) (v1) {$\zeta_0$};
    \node at (4,0) (cdots) {$\zeta_{t+1}$\ .};
    \draw (u0) edge node [above] {$t$} (v1);
    \draw (v1) edge node [above] {$t$} (cdots);
\end{tikzpicture}
\end{center}
}

\noindent
By Property 1, $m_0=2$, and by Property 2, $i=1$.

Conversely, suppose that condition $i)$ or $ii)$ holds. Then the corresponding path, as shown in Figure \ref{fig4}, allows us to apply Theorem \ref{tstcsrad} and conclude that $\rad(P(t))$ belongs to $\mathfrak{C}_{S(t)}$.

\begin{figure}[h!]
\begin{minipage}[c]{.5\textwidth}
\hspace{15mm} \begin{tikzpicture}[square/.style={regular polygon,regular polygon sides=4,minimum width=14mm,inner sep=0pt}]
\tikzset{vnode/.style={draw,thick,circle,minimum width=9mm,inner sep=0pt}};
    \node at (0,0) (u0) {$\zeta_{t+1}$};
    \node at (2,0) (v1) {$\zeta_0$};
    \node at (4,0) (cdots) {$\zeta_{t+1}$\ ,};
    \draw (u0) edge node [above] {$t$} (v1);
    \draw (v1) edge node [above] {$t$} (cdots);
\end{tikzpicture}
\end{minipage}\hspace{-28mm} 
\begin{minipage}[c]{.52\textwidth}
 \centering
 \begin{tikzpicture}[square/.style={regular polygon,regular polygon sides=4,minimum width=14mm,inner sep=0pt}]
\tikzset{vnode/.style={draw,thick,circle,minimum width=9mm,inner sep=0pt}};
    \node at (0,0) (u0) {$\zeta_0$};
    \node at (2,0) (v1) {$\zeta_{t+1}$};
    \node at (4,0) (cdots) {$\zeta_0$\ .};
    \draw (u0) edge node [above] {$t$} (v1);
    \draw (v1) edge node [above] {$t$} (cdots);
\end{tikzpicture}
\end{minipage}
\caption{The paths from conditions $i)$ and $ii)$, respectively.}\label{fig4}
\end{figure}
Hence, the connected component $\mathfrak{C}_{S(t)}$ of $\leftindex_s\Gamma_{\Lambda(n,\overline{m})}$ is $\Omega$-stable.
\end{proof}

\begin{remark}\label{rfdiij}\rm
\begin{enumerate}[$i)$]
    \item As a consequence of Lemma \ref{le}, all the components of $\leftindex_s\Gamma_{\Lambda(n,\overline{m})}$ containing non-periodic simple modules are $\Omega$-stable if and only if one of the following conditions holds.
\begin{enumerate}
    \item[$a)$] $\overline{m}=(2,\ldots,2)$;
    \item[$b)$] $i=1$, $m_0=2$ and $m_1>2$.
    \end{enumerate}
    \item Another consequence of Lemma \ref{le} and previous results is that the following classes of modules coincide for Brauer graph algebras associated to the star $\mathcal{W}_{n,\overline{m}}$.

\hspace{-24mm}
\begin{tabular}{p{1.5cm}lll}
    & \multirow{2}{*}{\parbox{50mm}{\{Periodic simple modules\} $=$}} & \hspace{-2mm}\multirow{2}{*}{\parbox{50mm}{\{Uniserial simple modules\} $=$}} & \hspace{-3mm}\multirow{2}{*}{\parbox{47mm}{\{Simple modules associated to exceptional edges\}}}\\
    & & \\
\end{tabular}
\end{enumerate}
\end{remark}

In the following subsections, we consider the cases provided in Lemma \ref{le} to calculate the universal deformation rings of modules in $\Omega$-stable components containing non-periodic simple modules.

\subsection{First case of $\Omega$-stability}

\begin{figure}[h!]
{\tiny
\begin{tikzcd}[column sep={30pt,between origins}, row sep=10pt]
	&&&&&&&& \ \rotatebox{18}{$\cdots$} \\
	&&&&&& \Omega\left(M\left[g_t\alpha_{t-1}\right]\right) \\
	\ \ \rotatebox{160}{$\cdots$} &&&& \ \Omega^{-1}(S(t)) \\
	&& S(t) \\
	\ \ \ \ \rotatebox{22}{$\cdots$} &&&& M\left[g_t\alpha_{t-1}\right]\ \\
	&&&&&& \ \ \ \rotatebox{162}{$\cdots$}\ \ 
	\arrow[from=2-7, to=1-9]
	\arrow[from=3-1, to=4-3]
	\arrow[from=3-5, to=2-7]
	\arrow[from=4-3, to=3-5]
	\arrow[from=4-3, to=5-5]
	\arrow[from=5-1, to=4-3]
	\arrow[from=5-5, to=6-7]
\end{tikzcd}
}
\caption{Component $\mathfrak{C}_{S(t)}$}\label{figch1}
\end{figure}
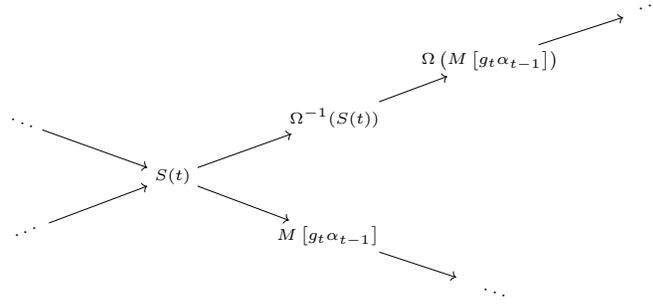

In this subsection, we consider the simple modules $S(t)$ in\hspace{-2mm} $\mod\hspace{-1.5mm} -\hspace{-.3mm} \Lambda(n,\overline{m})$ for which $m_{t+1}=2$. In Figure \ref{figch1}, we show some arrows in $\mathfrak{C}_{S(t)}$ to precise the action of $\Omega$ in $\mathfrak{C}_{S(t)}$, where
\vspace{-8mm}

\begin{center}
\begin{equation*}
g_t=\begin{cases}
e_n & \text{if $t=0$ and $i\neq n+1$,}\\
\delta_{t-1}^{-1} & \text{otherwise.}
\end{cases}
\end{equation*}
\end{center}

Hence, the modules in the diagonal
$$\xymatrix{\mathcal{D}_t: & \Omega^{-1}(S(t)) \ar[r] & \Omega^{-1}\left(M\left[\delta_{t-1}^{-1}\alpha_t\right]\right) \ar[r] & \cdots}$$
form a complete list of representatives of the $\Omega$-orbits in $\mathfrak{C}$. To explicitly describe the modules in the diagonal $\mathcal{D}_t$ we introduce some strings based on the $x_{l,p}$ from Def\mbox inition \ref{dsstr}.

\begin{defi}\label{dmorestrings} Consider $l,r\in\{0\ldots,n\}$. We def\mbox ine the strings $\mu_{l,r}$, $\mu_l$ and $\rho_{l,r}$ as follows.
\begin{enumerate}[$i)$]
    \item $\mu_{l,r}:=A_lx_{l,r}$;
    \item $\mu_l:=\mu_{l,l-1}$;
    \item $\rho_{l,r}:=\left(\delta_l^{-1}\mu_l\right) \cdots \left(\delta_1^{-1}\mu_1\right)\left(\delta_0^{-1}\mu_{0,r}\right)$ for $l \in \{0,...,i-1\}$;
    \item $\rho_l:=\rho_{l,l}$.
\end{enumerate}
\end{defi}

\begin{lemm}\label{lmjcomp} Suppose that $\overline{m}=(2,\ldots,2)$. Let $\mathfrak{C}$ be a connected component of $\leftindex_s\Gamma_{\Lambda(n,\overline{m})}$ containing a non-periodic simple module $S(t)$ and write
$$\xymatrix{\mathcal{D}_t: & M_0 \ar[r] & M_1 \ar[r] & \cdots}\ ,$$
where $\mathcal{D}_t$ is the diagonal given in Remark \ref{rfdiij}.$ii)$. Then, $M_j=M[C_j]$, where $C_j$ is a string of one of the following forms:
\begin{enumerate}
\item[$i)$] If $i\neq n+1$ and $0\leq j<t$, def\mbox ine $C_j=\left(\delta_t^{-1}\mu_t\right) \cdots \left(\delta_{t-j}^{-1}\mu_{t-j}\right)$.
\end{enumerate}    
\noindent    
For the cases ii) and iii) below, suppose that $j\geq t$ and write $j-t=q(n+1)+r$, where $q,r$ are non-negative integers and $r\leq n$.

\begin{enumerate}
\item[$ii)$] If $i\neq n+1$ and $0\leq r\leq n-i+1$, def\mbox ine
\begin{equation*}
C_j:= \begin{cases}
\rho_{t,n-r} & \text{if $q=0$,}\\
\rho_{t,i-1}\rho_{i-1,n-r} & \text{if $q=1$,}\\
\rho_{t,i-1}\rho_{i-1}^{q-1}\rho_{i-1,n-r} & \text{if $q\geq 2$.}
\end{cases}
\end{equation*}
\item[] If $i\neq n+1$ and $n-i+2\leq r\leq n$, def\mbox ine
\begin{equation*}
C_j:= \begin{cases}
\rho_{t,i-1} D & \text{if $q=0$,}\\
\rho_{t,i-1}\rho_{i-1}^q D & \text{if $q\geq 1$,}
\end{cases}
\end{equation*}
where $D:=\left(\delta_{i-1}^{-1}\mu_{i-1}\right) \cdots \left(\delta_{n-r+1}^{-1}\mu_{n-r+1}\right)$.
\item[$iii)$] If $i=n+1$, def\mbox ine
$$
C_j=(\delta_t^{-1}\mu_t)(\delta_{t-1}^{-1}\mu_{t-1})\cdots (\delta_{t-j}^{-1}\mu_{t-j}).
$$
\end{enumerate}
\end{lemm}

\begin{proof}
Using the injective envelope $S(t)\xhookrightarrow{} P(t)$, we obtain $\Omega^{-1}(S(t))=M\left[ \delta_t^{-1}\mu_t \right]=M\left[ C_0 \right]$. Similarly, from the injective envelope\mbox{} $M\left[ \delta_t^{-1}\mu_t \right]\xhookrightarrow{} P(t)\oplus P(t-1)$, we have
\begin{equation*}
\Omega^{-1}\left(M\left[g_t^{-1} \cdot \alpha_{t-1}\right]\right)=\begin{cases}

M\left[ \delta_0^{-1} \cdot A_0 \cdot x_{0,n-1}\right] & \text{if $t=0$  and $i\neq n+1$,}\\
M\left[ \delta_t^{-1} \cdot \mu_t \cdot \delta_{t-1}^{-1} \cdot \mu_{t-1} \right] & \text{if $t\neq 0$ or $i=n+1$.}
\end{cases}
\end{equation*}
\noindent
This implies that $\Omega^{-1}\left(M\left[g_t^{-1}\alpha_{t-1}\right]\right)=M\left[ C_1 \right]$. Then, $\mathcal{D}_t$ is the diagonal of $\mathfrak{C}$ starting with the arrow $M[C_0]\rightarrow M[C_1]$.

Now suppose that $i\neq n+1$. Observe that we have a sequence
\begin{equation}\label{eqb}
\xymatrix{\mathcal{E}: & M[C_0] \ar[r] & M[C_1] \ar[r] & M[C_2] \ar[r] & \cdots}
\end{equation}
of arrows in $\leftindex_s\Gamma_{\Lambda(n,\overline{m})}$ because the following identities hold for every $j\geq 1$:
\begin{enumerate}[$i)$]
    \item $C_j={}_h(C_{j-1})$, if $1\leq j\leq t$, or $r=0$, or $n-i+2\leq r\leq n$.
    \item $C_{j-1}={}_c(C_j)$, if $0<r\leq n-i+1$.
\end{enumerate}
For $j\geq 2$, write $\tau(M_j)=M[D_j]$. Note that $D_j$ is not of the form $\delta_t^{-1}E$ for any string $E$. Since all $C_j$ has this form, $\tau(M_j)\neq M_{j-2}$. Thus, $\mathcal{E}$ is a diagonal, i.e., $\mathcal{D}_t=\mathcal{E}$.

If $i=n+1$, we have a path $\mathcal{E}$ as in \eqref{eqb} satisfying $M[C_{j+1}]=\leftindex_h(M[C_j])$ for every $j>0$. Similarly $\tau(M_j)\neq M_{j-2}$ for every $j\geq 2$, implying $\mathcal{D}_t=\mathcal{E}$.
\end{proof}

\begin{theo}\label{taister} Suppose that $S(t)$ is a simple module in\hspace{-2mm} $\mod\hspace{-1.5mm} -\hspace{-.3mm} \Lambda(n,\overline{m})$ such that $m_{t+1}=2$. Then, a module $M$ in $\mathfrak{C}_{S(t)}$ has stable endomorphism ring isomorphic to $\Bbbk$ if and only if $M$ is in the $\Omega$-orbit of one of the modules $M_j$, with $0\leq j\leq n$.
\end{theo}

\begin{proof}
First, we will prove that $M_j$ has stable endomorphism ring isomorphic to $\Bbbk$ for $j\in \{0,\ldots n\}$. By the form of the strings $C_j$ in Lemma \ref{lmjcomp}, all the substrings of $C_j$ inducing a canonical endomorphism of $M[C_j]$ are paths.

Let $f$ be a canonical endomorphism of $M_j$ of the form $M_j\twoheadrightarrow M[p]\xhookrightarrow{} M_j$, where $p$ is a path. If $p=e_t$ and $f$ is given by the factorizations $C_j=\delta_t^{-1}e_tE$ and $C_j=e_t\delta_t^{-1}E$ for some string $E$, then there are canonical homomorphisms $\varphi:M_j\twoheadrightarrow M[\delta_t^{-1}]\xhookrightarrow{} P(t)$ and $\psi:P(t)\twoheadrightarrow P(t)/\soc(P(t))\xhookrightarrow{} M_j$ such that $f=\psi\varphi$. Hence $f\in \mathcal{P}_{\Lambda(n,\overline{m})}(M_j,M_j)$. If $f$ is not given by such factorizations, then $f$ is given by factorizations of the form $C_j=F_1pp'F_2$ and $C_j=F_3p''pF_4$, where $F_1$, $F_2$, $F_3$ and $F_4$ are strings, and $p'$ and $p''$ are paths, such that $F_1$ and $F_3$ are not of the form $F_1'\beta_1$ for any string $F_1'$ and $\beta_1 \in Q_1$, and $F_2$ and $F_4$ are not of the form $\beta_2 F_2'$ for any string $F_2'$ and $\beta_2 \in Q_1$. Then $f=\psi\varphi-f'$, where $\varphi$ and $\psi$ are canonical homomorphisms of the form $M_j\twoheadrightarrow M[pp'w]\xhookrightarrow{} P(t(p'))$ and $P(t(p'))\twoheadrightarrow M[w'up]\xhookrightarrow{} M_j$, respectively, with
\begin{equation*}
w:= \begin{cases}
e_{t(p')} & \text{if $F_2$ is a trivial string,}\\
\delta_{t(p')}^{-1} & \text{otherwise,}
\end{cases}
\end{equation*}
$u$ is the path such that $A_{t(p')}^2=upp'$, the string $w'$ is given by
\begin{equation*}
w':= \begin{cases}
e_{t(p')} & \text{if $u\neq p''$,}\\
\delta_{t(p')}^{-1} & \text{if $u=p''$,}
\end{cases}
\end{equation*}
and $f'$ is the zero homomorphism if $u\neq p''$ (in particular, it happens when $p$ is a trivial path) and $f$ is the canonical homomorphism $M_j\twoheadrightarrow S(t(p'))\xhookrightarrow{} M_j$ if $u=p''$. Hence, every canonical endomorphism of $M_j$ dif\mbox ferent from the identity is in $\mathcal{P}(M_j,M_j)$. Consequently, we get $\send_{\Lambda(n,\overline{m})}(M_j)\cong \Bbbk$ for $0\leq j\leq n$.

Now, we will prove that $\send_{\Lambda(n,\overline{m})}(M_j)\not\cong \Bbbk$ for $j>n$. Consider the canonical homomorphism $g:M_j\twoheadrightarrow S(t)\xhookrightarrow{} M_j$ induced by the factorizations $C_j=(C_n\delta_t^{-1})e_tG$ and $C_j=e_tC_j$, where $G$ is a string. Suppose that $g\in \mathcal{P}(M_j,M_j)$. Then, there are a projective indecomposable module $P$ and homomorphisms $\varphi: M_j\rightarrow P$ and $\psi: P\rightarrow M_j$ such that $S(t)$ is a submodule of $Im(\psi\varphi)$. By the form of $C_j$, we have that $P\cong P(t)$ and $\psi$ is generated by the canonical homomorphism $P(t)\twoheadrightarrow M[\delta_t^{-1}\mu_t]\xhookrightarrow{} M_j$, and by the forms of $C_j$ and $P(t)$, we have that $\varphi$ is generated by the canonical homomorphism $M_j\twoheadrightarrow M[\mu_{t+1}\delta_t^{-1}]\xhookrightarrow{} P(t)$. Because the composition of these two homomorphisms sends the top submodule $S(t)\oplus M[\mu_{t+1,t-1}]$ of $M_j$ to a submodule of $M_j$ of the same form, we deduce that $M[\mu_{t+1,t-1}]$ is not contained in $\ker(g)$, which is a contradiction with the choice of $g$. Therefore, $g\notin \mathcal{P}(M_j,M_j)$ and hence $\send_{\Lambda(n,\overline{m})}(M_j)\not\cong \Bbbk$ for $j>n$.

Finally, the statement follows from Remark \ref{rfdiij} part $ii)$.
\end{proof}

\begin{coro} Suppose that $S(t)$ is a simple module in\hspace{-2mm} $\mod\hspace{-1.5mm} -\hspace{-.3mm} \Lambda(n,\overline{m})$ such that $m_{t+1}=2$. Then, there are exactly $2n+2$ $\Omega^2$-orbits in $\mathfrak{C}$ containing modules with stable endomorphism ring isomorphic to $\Bbbk$.
\end{coro}

\begin{theo}\label{tex2222} Suppose that $S(t)$ is a simple module in\hspace{-2mm} $\mod\hspace{-1.5mm} -\hspace{-.3mm} \Lambda(n,\overline{m})$ such that $m_{t+1}=2$. If $M$ is a module in $\mathfrak{C}_{S(t)}$ with stable endomorphism ring isomorphic to $\Bbbk$, then
\begin{equation*}
\dim_{\Bbbk} \Ext_{\Lambda(n,\overline{m})}(M,M)=
\begin{cases}
1 & \text{ if $M$ is in the same $\Omega$-orbit as $S(t)$ or $M_n$,}\\
0 & \text{ otherwise.}
\end{cases}
\end{equation*}
\end{theo}

\begin{proof}
First, we will prove that the statement is true for the modules $M_j$ in the diagonal $\mathcal{D}_t$.

For $j=0$, $\Omega(M_0)=S(t)$. The natural inclusion $S(t)\xhookrightarrow{} M_0$ generates $\text{Hom}_{\Lambda(n,\overline{m})}(\Omega(M_0),M_0)$ and does not factor through any projective module. Hence,
$$
\dim_{\Bbbk} \Ext_{\Lambda(n,\overline{m})}^1(M_0,M_0)=\dim_{\Bbbk} \shom_{\Lambda(n,\overline{m})}(\Omega(M_0),M_0)=1.
$$

Now, suppose that $j\neq 0$. For each $a\in \{0,\ldots,i-1\}$ and for each $b\in \{a,a+1,\ldots,i-1\}$ we set
$$\gamma_{a,b}:=(\delta_a^{-1}\alpha_a)(\delta_{a+1}^{-1}\alpha_{a+1})\cdots(\delta_b^{-1}\alpha_b).$$
Taking a projective cover $P$, we obtain that $\Omega(M_j)=M[D_j]$, where
\vspace{-4mm}

\begin{equation*}
D_j:=\begin{cases}
\gamma_{t-j,t-1} & \text{ if $0 \leq j\leq t$,}\\
x_{n-j+t+1,0}\gamma_{0,t-1} & \text{ if $t+1 \leq j\leq n-i+t+1$,}\\
\delta_{i-1}^{-1}x_{i-1,0}\gamma_{0,t-1} & \text{ if $j=n-i+t+2$,}\\
\gamma_{n-j+t+1,i-2}\delta_{i-1}^{-1}x_{i-1,0}\gamma_{0,t-1} & \text{ if $n-i+t+3 \leq j\leq n$.}
\end{cases}
\end{equation*}

If $j\neq n$, then all the canonical homomorphisms from $\Omega(M_j)$ to $M_j$ are of the form $\Omega(M_j)\twoheadrightarrow S(k)\xhookrightarrow{} M_j$, where $k\in \{0,\ldots,i-1\}$. Moreover, there is at most one homomorphism of this form for each value of $k$, which we denote by $f_k$. Let $k_1,\ldots, k_s$ be all the vertices in $Q_0$ that induce canonical homomorphisms of this form. For every $l\in \{1,\ldots,s\}$, we have that $D_j=E_1e_{k_l}p_lE_2$, for some strings $E_1$ and $E_2$, and some path $p_l$, such that $E_1$ is not of the form $E'\beta$ for any string $E'$ and for any arrow $\beta$, and $E_2$ is given by
\begin{equation*}
E_2:= \begin{cases}
e_{t-1} & \text{if $l=s$,}\\
\delta_{t(p_l)}^{-1}e_{k_{l+1}}\alpha_{k_{l+1}}E_2' & \text{for a string $E_2'$, if $l<s$.}
\end{cases}
\end{equation*}
Under these assumptions, $f_{k_l}$ is induced by factorizations $D_j=E_1e_{k_l}p_lE_2$ and $C_j=E_3p'e_{k_l}E_4$ for some strings $E_3$ and $E_4$, and some path $p'$, such that $E_3$ is not of the form $E'\beta_1$ for any string $E'$ and for any arrow $\beta_1$, and $E_4$ is not of the form $\beta_2E''$ for any string $E''$ and for any arrow $\beta_2$. We denote by $\varphi_{k_l}$ and $\psi_{k_l}$ the canonical homomorphisms $\Omega(M_j)\twoheadrightarrow M[p_lw]\xhookrightarrow{} P(t(p_l))$ and $P(t(p_l))\twoheadrightarrow M[\delta_{t(p_l)}^{-1}p']\xhookrightarrow{} M_j$, respectively, where
\begin{equation*}
w:= \begin{cases}
e_{t-1} & \text{if $l=s$,}\\
\delta_{t(p_l)}^{-1} & \text{if $l<s$.}
\end{cases}
\end{equation*}
Hence, $f_{k_s}=\psi_{k_s}\varphi_{k_s}$ and $f_{k_l}=\psi_{k_l}\varphi_{k_l}-f_{k_{l+1}}$ for $l<s$, which implies that all the canonical homomorphisms from $\Omega(M_j)$ to $M_j$ are in $\mathcal{P}(\Omega(M_j),M_j)$. Therefore,
$$
\dim_{\Bbbk} \Ext_{\Lambda(n,\overline{m})}^1(M_j,M_j)=\dim_{\Bbbk} \shom_{\Lambda(n,\overline{m})}(\Omega(M_j),M_j)=0
$$
for $1\leq j<n$.

For $j=n$, we consider the canonical homomorphism $f:\Omega(M_n)\twoheadrightarrow M[\alpha_{t-1}]\xhookrightarrow{} M_n$ induced by factorizations $D_j=F_1\alpha_{t-1}$ and $C_j=F_2\alpha_{t-1}$, for some strings $F_1$ and $F_2$. For every pair of canonical homomorphisms $\varphi:M_n\rightarrow P$ and $\psi:P\rightarrow M_n$, with $P$ a projective indecomposable module, such that $M[\alpha_{t-1}]\subseteq \im(\varphi)$, we obtain that $P\cong P(t)$. In this case, we have that $\im(\varphi)\cap \soc(P(t))\neq 0$ and hence $\im(\varphi)\cap \ker(\psi)\neq 0$. This implies that $\im(f)\cap \ker(\psi\varphi)\neq 0$ and therefore $f\notin \mathcal{P}(\Omega(M_n),M_n)$. The other canonical homomorphisms from $\Omega(M_n)$ to $M_n$ are of the form $\Omega(M_n)\twoheadrightarrow S(k)\xhookrightarrow{} M_n$, where $k\in \{0,\ldots,i-1\}$ and the same argument as for the case $j<n$ shows that all these homomorphisms are in $\mathcal{P}(\Omega(M_n),M_n)$. This proves that
$$
\dim_{\Bbbk} \Ext_{\Lambda(n,\overline{m})}^1(M_j,M_j)=\dim_{\Bbbk} \shom_{\Lambda(n,\overline{m})}(\Omega(M_j),M_j)=1.
$$

Finally, the statement follows for all the modules in $\mathfrak{C}$ because the modules in $\mathcal{D}_t$ form a complete list of representatives of the $\Omega$-orbits in $\mathfrak{C}$ and $\Ext_{\Lambda(n,\overline{m})}^1(M,M)\cong \Ext_{\Lambda(n,\overline{m})}^1(\Omega(M),\Omega(M))$ for every $M$ in\hspace{-2mm} $\mod$\mbox- $\Lambda(n,\overline{m})$.

\end{proof}

\begin{theo}\label{udrs1} Suppose that $S(t)$ is a simple module in\hspace{-2mm} $\mod\hspace{-1.5mm} -\hspace{-.3mm} \Lambda(n,\overline{m})$ such that $m_{t+1}=2$. If $M$ is a module in $\mathfrak{C}_{S(t)}$ with stable endomorphism ring isomorphic to $\Bbbk$, then
\vspace{-5mm}

\begin{center}
\begin{equation*}
R(\Lambda(n,\overline{m}),M)=
\begin{cases}
\Bbbk\llbracket x\rrbracket/\langle x^2\rangle & \text{ if $M$ is in the same $\Omega$-orbit as $S(t)$,}\\
\Bbbk\llbracket x\rrbracket & \text{ if $M$ is in the same $\Omega$-orbit as $M_n$,}\\
\Bbbk & \text{ otherwise.}
\end{cases}
\end{equation*}
\end{center}
\end{theo}

\begin{proof}
By Theorem \ref{trsyz} and Remark \ref{rfdiij}, it suf\mbox f\mbox ices to consider the modules $M_j$ in $\mathcal{D}_t$.

For $0<j<n$, Theorem \ref{tex2222} and \cite[Proposition 2.1]{bv1} imply that $R(\Lambda(n,\overline{m}),M_j)\cong \Bbbk$.

For $j=0$, by Theorem \ref{trsyz}, we can consider $\Omega(M_0)=S(t)$ (instead of $M_0$). From Theorem \ref{tex2222} and \cite[Proposition 2.1]{bv1}, we know that $R(\Lambda(n,\overline{m}),S(t))$ is a quotient of the power series ring $\Bbbk\llbracket x\rrbracket$. Consider the sequence $\mathcal{L}=\{V_0,V_1\}$, where $V_0:=S(t)$ and $V_1:=M[\delta_t^{-1}]$ joint with the natural inclusion $\iota:V_0\xhookrightarrow{} V_1$ and the canonical projection $\epsilon:V_1\twoheadrightarrow V_0$. Set $\sigma:=\iota\epsilon$. Then, $\ker(\sigma)\cong V_0$ and $\im(\sigma)\cong V_0$. Moreover, $\Hom_{\Lambda(n,\overline{m})}(V_1,V_0)$ has dimension $1$ as vector space because it is generated by $\epsilon$, and
$$\Ext_{\Lambda(n,\overline{m})}(V_1,V_0)\cong \shom_{\Lambda(n,\overline{m})}(\Omega(V_1),V_0)=\shom_{\Lambda(n,\overline{m})}(M[\mu_{t+1}],S(t))=0$$
because there are no canonical homomorphisms from $M[\mu_{t+1}]$ to $S(t)$. Hence, by Theorem \ref{anote}, we obtain $R(\Lambda(n,\overline{m}),M_j)\cong \Bbbk\llbracket x\rrbracket/\langle x^2\rangle$.

For $j=n$, we can work with $\Omega(M_n)=M[D_n]$. We def\mbox ine the sequence $\mathcal{L}':=\{W_0,W_1,\ldots\}$ by taking $W_0=M[D_n]$ and $W_j=M[D_n(\alpha_tD_n)^j]$ for $j>0$. Also, let $\iota_j:W_{j-1}\xhookrightarrow{} W_j$ be a natural inclusion and $\epsilon_j:W_j\twoheadrightarrow W_{j-1}$ be a canonical projection, and put $\sigma_j:=\iota_j\epsilon_j$. Then, it holds that $\ker{\sigma_j}\cong  W_0$ and $\im(\sigma_j^j)\cong W_0$. Since $\mathcal{L}'$ is inf\mbox inite, Theorem \ref{anote} implies that
$$R\left(\Lambda(n,\overline{m}),M_n\right)\cong R\left(\Lambda(n,\overline{m}),M[D_n]\right)\cong \Bbbk\llbracket x\rrbracket.$$

Finally, the statement follows from Theorem \ref{trsyz} part $ii)$.
\end{proof}

\subsection{Second case of $\Omega$-stability}

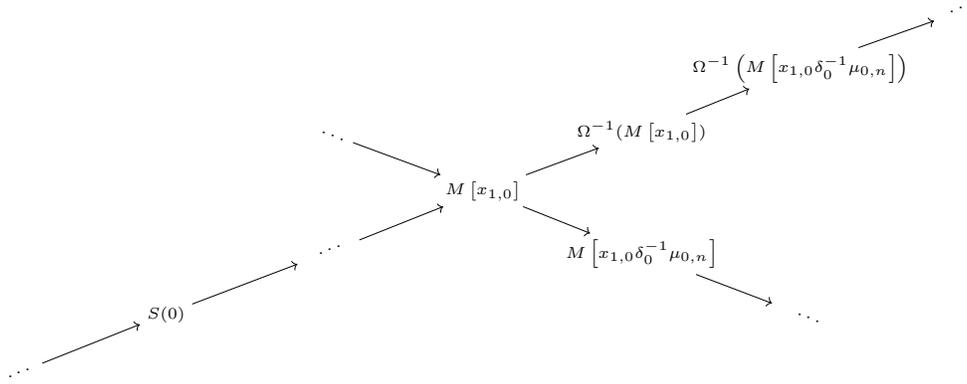
\begin{figure}[h!]
{\tiny
\begin{tikzcd}[column sep={30pt,between origins}, row sep=10pt]
	&&&&&&&&&&&& \ \rotatebox{22}{$\cdots$} \\
	&&&&&&&& {} && {\Omega^{-1}\left(M\left[x_{1,0}\delta_0^{-1}\mu_{0,n}\right]\right)} \\
	&&&& \ \ \rotatebox{161}{$\cdots$} &&&& {\Omega^{-1}(M\left[x_{1,0}\right])} \\
	&&&&&& {M\left[x_{1,0}\right]} \\
	&&&& \ \ \ \rotatebox{21}{$\cdots$\ } &&&& {M\left[x_{1,0}\delta_0^{-1}\mu_{0,n}\right]} \\
	&& {S(0)} &&&&&&&& \ \ \rotatebox{160}{$\cdots$} \\
	\ \ \ \ \rotatebox{23}{$\cdots$}
	\arrow[from=2-11, to=1-13]
	\arrow[from=3-5, to=4-7]
	\arrow[from=3-9, to=2-11]
	\arrow[from=4-7, to=3-9]
	\arrow[from=4-7, to=5-9]
	\arrow[from=5-5, to=4-7]
	\arrow[from=5-9, to=6-11]
	\arrow[from=6-3, to=5-5]
	\arrow[from=7-1, to=6-3]
\end{tikzcd}
}
\caption{Component $\mathfrak{C}_{S(0)}$}\label{figch2}
\end{figure}

In this subsection, we consider the star $\mathcal{W}_{n,\overline{m}}$ and we assume that $i=1$, $m_0=2$ and $m_1>2$. In this case, the unique non-periodic simple module over $\Lambda(n,\overline{m})$ is $S(0)$. In Figure \ref{figch2} we show some arrows of the connected component of $\leftindex_s\Gamma_{\Lambda(n,\overline{m})}$ corresponding to $S(0)$.

By the action of the operator $\Omega$, the modules in the diagonal
$$\xymatrix{\mathcal{D}_0: & M\left[x_{1,0}\right] \ar[r] & M\left[x_{1,0}\delta_0^{-1}\mu_{0,n}\right] \ar[r] & \cdots}$$
form a complete list of representatives of the $\Omega$-orbits in $\mathfrak{C}$. Observe that this diagonal can be written as
$$\xymatrix{\mathcal{D}_0: & M_0 \ar[r] & M_1 \ar[r] & \cdots},$$
with $M_j=M[C_j]$ for $j\in \{0,1,\ldots\}$, where the strings $C_j$ are def\mbox ined as follows:
\begin{enumerate}[$i)$]
    \item $C_0:=x_{1,0}$.
    \item For $j\geq 1$, write $j-1=q(n+1)+r$, with $q,r\in \mathbb{Z}$ and $0\leq r\leq n$, and def\mbox ine
    \begin{equation*}
C_j:= \begin{cases}
x_{1,0}\delta_0^{-1}\mu_{0,n-j+1} & \text{if $q=0$,}\\
x_{1,0}(\delta_0^{-1}A_0)^q\mu_{0,n-r+1} & \text{if $q>0$.}
\end{cases}
\end{equation*}
Note that $C_j={}_h(C_{j-1})$ if $r=0$, and $C_{j-1}={}_c(C_j)$ if $r\neq 0$.
\end{enumerate}

\begin{theo}\label{ti1m02m12}
Suppose that $i=1$, $m_0=2$ and $m_1>2$. Then, a module $M$ in $\mathfrak{C}_{S(0)}$ has stable endomorphism ring isomorphic to $\Bbbk$ if and only if $M$ is in the $\Omega$-orbit of one of the modules $M_j$, with $0\leq j\leq n$.
\end{theo}

\begin{proof}
For $j=0$, we have that $\send_{\Lambda(n,\overline{m})}(M_0,M_0)$ is generated by the identity and hence we have that $\send_{\Lambda(n,\overline{m})}(M_0,M_0)\cong \Bbbk$.

Assume that $1\leq j\leq n$. Then all the canonical homomorphisms from $M_j$ to itself are of the form $M_j\twoheadrightarrow S(0)\xhookrightarrow{} M_j$, $M_j\twoheadrightarrow M[x_{1,n-j+1}]\xhookrightarrow{} M_j$ and $M_j\twoheadrightarrow M[x_{0,n-j+1}]\xhookrightarrow{} M_j$. The canonical homomorphism $f:M_j\twoheadrightarrow S(0)\xhookrightarrow{} M_j$ factors as $f=\psi \varphi$, where $\varphi$ and $\psi$ are the canonical homomorphisms $M_j\twoheadrightarrow M[\mu_{0,n-j+1}]\xhookrightarrow{} P(n-j+1)$ and $P(n-j+1)\twoheadrightarrow M[x_{n-j+1,0}]\xhookrightarrow{} M_j$, respectively; the canonical homomorphism $g:M_j\twoheadrightarrow M[x_{1,n-j+1}]\xhookrightarrow{} M_j$ can be expressed as $g=\psi'\varphi'-f$, where $\varphi'$ and $\psi'$ are the canonical homomorphisms $M_j\twoheadrightarrow M[x_{1,0}\delta_0^{-1}]\xhookrightarrow{} P(0)$ and $P(0)\twoheadrightarrow M[\delta_0^{-(m_1-1)}\mu_{0,n-j+1}]\xhookrightarrow{} M_j$, respectively; and the canonical homomorphism $h:M_j\twoheadrightarrow M[x_{0,n-j+1}]\xhookrightarrow{} M_j$ factors as $g=\psi''\varphi''$, where $\varphi''$ and $\psi''$ are the canonical homomorphisms $M_j\twoheadrightarrow M[\mu_{0,n-j+1}]\xhookrightarrow{} P(n-j+1)$ and $P(n-j+1)\twoheadrightarrow M[A_{n-j+1}]\xhookrightarrow{} M_j$, respectively. Since all the canonical homomorphisms from $M_j$ to $M_j$ distinct from the identity are in $\mathcal{P}(M_j,M_j)$, we obtain that $\send_{\Lambda(n,\overline{m})}(M_j,M_j)\cong \Bbbk$ for $1\leq j\leq n$.

For $j>n$, the canonical homomorphism $M_j\twoheadrightarrow S(0)\xhookrightarrow{} M_j$ induced by factorizations of the form $C_j=x_{1,0}\delta_0^{-1}e_0E$ and $C_j=x_{1,0}e_0E'$, with $E$ and $E'$ strings, does not belong to $\mathcal{P}(M_j,M_j)$ due to the form of the projective indecomposable modules.
\end{proof}

\begin{theo}\label{texti1m12}
Suppose that $i=1$, $m_0=2$ and $m_1>2$. If $M$ is a module in $\mathfrak{C}_{S(0)}$ with stable endomorphism ring isomorphic to $\Bbbk$, then
\begin{equation*}
\dim_{\Bbbk} \Ext_{\Lambda(n,\overline{m})}(M,M)=
\begin{cases}
1 & \text{ if $M$ is in the same $\Omega$-orbit as $M_0$ or $M_n$,}\\
0 & \text{ otherwise.}
\end{cases}
\end{equation*}
\end{theo}

\begin{proof}
For $j=0$, the vector space
$$\Ext_{\Lambda(n,\overline{m})}(M_0,M_0)\cong \shom_{\Lambda(n,\overline{m})}(\Omega(M_0),M_0)\cong \shom_{\Lambda(n,\overline{m})}(M[A_1],M[x_{1,0}])$$
is generated by a canonical projection which does not belong to $\mathcal{P}(\Omega(M_0),M_0)$ because of the length of the radical series of $P(1)$. Therefore, $\dim_{\Bbbk} \Ext_{\Lambda(n,\overline{m})}(M_0,M_0)=1$.

For $1\leq j\leq n-1$, the canonical homomorphisms from $\Omega(M_j)=M[x_{n-j+2,0}\delta_0^{-(m_1-1)}\mu_{0,1}]$ to $M_j$ are $f:\Omega(M_j)\twoheadrightarrow S(0)\xhookrightarrow{} M_j$ and $g:\Omega(M_j)\twoheadrightarrow M[x_{0,n-j+1}]\xhookrightarrow{} M_j$ for which we have that $f=\psi \varphi$ and $g=\psi' \varphi$, where $\varphi$, $\psi$ and $\psi'$ are the following canonical homomorphisms
$$\varphi:\Omega(M_j)\twoheadrightarrow M[\mu_{0,1}]\xhookrightarrow{} P(1),$$
$$\psi:P(1)\twoheadrightarrow M[x_{1,0}]\xhookrightarrow{} M_j,$$
$$\psi':P(1)\twoheadrightarrow M[\mu_{1,n-j+1}]\xhookrightarrow{} M_j.$$
Hence, $\dim_{\Bbbk} \Ext_{\Lambda(n,\overline{m})}(M_j,M_j)=0$ for $1\leq j\leq n-1$.

For $j=n$, the canonical homomorphisms from $\Omega(M_n)=M[x_{2,0}\delta_0^{-(m_1-1)}\mu_{0,1}]$ to $M_n$ are $f':\Omega(M_n)\twoheadrightarrow S(0)\xhookrightarrow{} M_n$, $g':\Omega(M_n)\twoheadrightarrow M[\alpha_0]\xhookrightarrow{} M_j$ and $h:\Omega(M_j)\twoheadrightarrow M[\delta_0^{-1}\mu_{0,1}]\xhookrightarrow{} M_j$. The same argument used in the case $1\leq j\leq n-1$ shows that $f',g'\in \mathcal{P}(\Omega(M_n),M_n)$, whereas we observe that $h\not\in \mathcal{P}(\Omega(M_n),M_n)$ due to the form of the projective indecomposable modules. This implies that $\dim_{\Bbbk} \Ext_{\Lambda(n,\overline{m})}(M_n,M_n)=1$.
\end{proof}

\begin{theo}\label{udrsc2} Suppose that $i=1$, $m_0=2$ and $m_1>2$. If $M$ is a module in $\mathfrak{C}_{S(0)}$ with stable endomorphism ring isomorphic to $\Bbbk$, then
\vspace{-5mm}

\begin{center}
\begin{equation*}
R(\Lambda(n,\overline{m}),M)=
\begin{cases}
\Bbbk\llbracket x\rrbracket/\langle x^2\rangle & \text{ if $M$ is in the same $\Omega$-orbit as $M_0$,}\\
\Bbbk\llbracket x\rrbracket/\langle x^{m_1}\rangle & \text{ if $M$ is in the same $\Omega$-orbit as $S(0)$,}\\
\Bbbk & \text{ otherwise.}
\end{cases}
\end{equation*}
\end{center}
\end{theo}

\begin{proof}
The proof of the statement for $M_j$, with $0\leq j\leq n-1$ is analogous to that of Theorem \ref{udrs1}. It is worth to say that for $j=0$ we can work directly with $M_0$ and we can take $V_0:=M[x_{1,0}]$, $V_1:=M[\mu_{1,0}]$ for which we have $\Omega(V_1)=S(1)$.

For $j=n$, we can work with $S(0)$ instead of $M_n$ since they are in the same $\Omega$-orbit. Let $\mathcal{L}:=\{W_0,\ldots,W_{m_1-1}\}$, where $W_0:=S(0)$ and $W_j:=M[\delta_0^{-j}]$ for $1\leq j\leq m_1-1$. Also, for $1\leq j\leq m_1-1$, let $\iota_j:W_{j-1}\xhookrightarrow{} W_j$ be the natural inclusion and and $\epsilon_j:W_j\twoheadrightarrow W_{j-1}$ be the canonical projection, and put $\sigma_j:=\iota_j\epsilon_j$. Then, it holds that $\ker{\sigma_j}\cong  W_0$ and $\im(\sigma_j^j)\cong W_0$. Moreover, $\dim_{\Bbbk} \Hom_{\Lambda(n,\overline{m})}(W_{m_1-1},W_0)=1$ because the unique canonical homomorphism from $W_{m_1-1}$ to $W_0$ is the canonical projection, and
$$\Ext_{\Lambda(n,\overline{m})}^1(W_{m_1-1},W_0)\cong \shom_{\Lambda(n,\overline{m})}(\Omega(W_{m_1-1}),W_0)=\shom_{\Lambda(n,\overline{m})}(M[\mu_{1,0}],S(0))=0$$
because there are no canonical homomorphisms from $M[\mu_{1,0}]$ to $S(0)$. In consequence, by Theorem \ref{anote}, $R(\Lambda(n,\overline{m}),S(0))=\Bbbk\llbracket x\rrbracket/\langle x^{m_1}\rangle$.

Finally, the statement follows from Theorem \ref{trsyz} part $ii)$.
\end{proof}

Applying Theorem \ref{tderst}, Theorem \ref{tedstm}, Lemma \ref{lbdes}, Theorem \ref{taister}, Theorem \ref{udrs1}, Theorem \ref{texti1m12} and Theorem \ref{udrsc2}, we obtain the following result for any generalized Brauer tree algebra.

\begin{coro}\label{czaingbt} Let $\mathcal{G}=(V,E,\mathfrak{m},\mathfrak{o})$ be a generalized Brauer tree and $\Lambda_{\mathcal{G}}$ be its corresponding Brauer graph algebra. Suppose that $\Lambda_{\mathcal{G}}$ is of non-polynomial growth. Def\mbox ine $l$ to be the number of vertices in $V$ with multiplicity $2$ and set $n:=|E|-1$.
\begin{enumerate}[$i)$]
    \item If $l>1$, then there exists a list $\{\mathfrak{C}_1,\ldots,\mathfrak{C}_{l-1}\}$ of $\Omega$-stable components of type $\mathbb{Z}A_{\infty}^{\infty}$ in $\leftindex_s\Gamma_{\Lambda_{\mathcal{G}}}$ such that, for every $i\in\{1,\ldots,l-1\}$, there are exactly $n+1$ $\Omega$-orbits of modules in $\mathfrak{C}_i$ with stable endomorphism ring isomorphic to $\Bbbk$. Moreover, a complete family of representative modules of these $\Omega$-orbits are given by a sectional path of the form\\
    $$\xymatrix{\mathcal{D}_t: & M_0 \ar[r] & M_1 \ar[r] & \cdots \ar[r] & M_n},$$\\
    where
\vspace{-8mm}

\begin{center}
\begin{equation*}
R(\Lambda_{\mathcal{G}},M_j)=
\begin{cases}
\Bbbk\llbracket x\rrbracket/\langle x^2\rangle & \text{ if $j=0$,}\\
\Bbbk\llbracket x\rrbracket & \text{ if $j=n$,}\\
\Bbbk & \text{ otherwise.}
\end{cases}
\end{equation*}
\end{center}

    \item If $l=1$ and there is a unique vertex $v\in V$ of multiplicity grater than $2$, then there is a $\Omega$-stable component $\mathfrak{C}$ of type $\mathbb{Z}A_{\infty}^{\infty}$ in $\leftindex_s\Gamma_{\Lambda_{\mathcal{G}}}$ such that there are exactly $n+1$ $\Omega$-orbits of modules in $\mathfrak{C}$ with stable endomorphism ring isomorphic to $\Bbbk$. Moreover, a complete family of representative modules of the $\Omega$-orbits are given by a sectional path of the form\\
    $$\xymatrix{\mathcal{D}_t: & M_0 \ar[r] & M_1 \ar[r] & \cdots \ar[r] & M_n},$$\\
    where

\vspace{-8mm}
\begin{center}
\begin{equation*}
R(\Lambda_{\mathcal{G}},M_j)=
\begin{cases}
\Bbbk\llbracket x\rrbracket/\langle x^2\rangle & \text{ if $j=0$,}\\
\Bbbk\llbracket x\rrbracket/\langle x^{\mathfrak{m}(v)}\rangle & \text{ if $j=n$,}\\
\Bbbk & \text{ otherwise.}
\end{cases}
\end{equation*}
\end{center}

\end{enumerate}
\end{coro}

\section*{Acknowledgments}

The authors thank Yadira Valdivieso-Díaz (UDLAP, Puebla, México) for valuable discussions related to this work.

\section*{Ethical Approval}

Not applicable

\section*{Funding}

The f\mbox irst author was partially supported by CODI (Universidad de Antioquia, UdeA), by project number 2020-33305. The second author was partially supported by CODI (Universidad de Antioquia, UdeA), by project numbers 2022-52654 and 2023-62291. The third author was supported by the research group PROMETE-KONRAD (Project No. 5INV1232) in Facultad de Matem\'aticas e Ingenier\'{\i}as at the Fundaci\'on Universitaria Konrad Lorenz, Bogot\'a, Colombia, by the Faculty Scholarship of the Of\mbox f\mbox ice of Academic Af\mbox fairs at the Valdosta State University, GA, USA, and by the research group \'Algebra UdeA in the Instituto de Matem\'aticas  at the Universidad de Antioquia in Medell\'{\i}n, Colombia. 

\section*{Availability of data and materials}  

Not applicable.


\begin{thebibliography}{YYYY}

\bibitem{az} Antipov, M.; Zvonareva A. \emph{Brauer graph algebras are closed under derived equivalence.} Math. Z. (2022). DOI: https://doi.org /10.1007/s00209-021-02937-x.

\bibitem{ass} Assem, I.; Simpson, D.; Skowro\'nski, A. \emph{Elements of the representation theory of associative algebras}. London Math. Soc. Student Texts 65, Cambridge University Press, Cambridge. (2006).

\bibitem{bc5} Bleher, F.; Chinburg, T. \emph{Deformations and derived categories}. Ann. Inst. Fourier (Grenoble) 55, 2285–2359. (2005).

\bibitem{bv1} Bleher, F. M.; Vélez-Marulanda, J. A. \emph{Universal deformation rings of modules over Frobenius algebras}. J. Algebra {\bf 367}, 176-202. (2012).

\bibitem{bv2} Bleher, F. M.; Vélez-Marulanda, J. A. \emph{Deformations of complexes for f\mbox inite-dimensional algebras}. J. Algebra 491, 90–140. (2017).

\bibitem{bw} Bleher, F. M.; Wackwitz, D. J. \emph{Universal deformation rings and self-injective Nakayama algebras}. J. Pure Appl. Algebra, 223 (1), pp. 218-244. (2019).

\bibitem{bs} Bocian, R.; Skowro\'nski, A. \emph{Symmetric special biserial algebras of Euclidean type}. Colloq. Math. 96(1), 121–148. (2003).

\bibitem{brr} Brou\'e, M. \emph{Equivalences of blocks of group algebras}. In: Finite Dimensional Algebras and Related Topics, Ottawa, ON, 1992, in: NATO Adv. Sci. Inst. Ser. C Math. Phys. Sci., vol.424, Kluwer Acad. Publ., Dordrecht, pp.1–26. (1994).

\bibitem{br} Butler, M. C. R.; Ringel, C. M. \emph{Auslander-Reiten sequences with a few middle terms and applications to string algebras}. Comm. Algebra {\bf 15}, 145-179. (1987).

\bibitem{cgrv} Calderón-Henao, Y.; Giraldo, H.; Rueda-Robayo, R.; Vélez-Marulanda, J. A. \emph{Universal deformation rings of string modules over certain class of self-injective special biserial algebras}. Comm. Algebra 47(12):4969–4988. (2019). DOI: 10.1080/00927872.2019.1609010.

\bibitem{duff} Duf\mbox f\mbox ield, D. \emph{Auslander-Reiten components of symmetric special biserial algebras}. J. Algebra. 508: 475–511. (2018). DOI: 10.1016/j.jalgebra.2018.03.040. 

\bibitem{kebt} Erdmann, K. \emph{Blocks of tame representation type and related algebras}, Springer Lecture Notes in Mathematics, {\bf 1428}. (1990).

\bibitem{kean1} Erdmann, K. \emph{Algebras with non-periodic bounded modules}. J. Algebra 475, p 308-326. (2017).

\bibitem{kes} Erdmann, K.; Skowro\'nski A. \emph{On Auslander–Reiten components of blocks and self-injective biserial algebras}. Trans. Amer. Math. Soc. 330(1), 165–189. (1992).

\bibitem{adr} Fonce-Camacho, A.; Giraldo, H.; Rizzo, P.; Vélez-Marulanda, J. A. \emph{On a deformation theory of f\mbox inite dimensional modules over repetitive algebras}. Algebr Represent Theor 26, 1–22. (2023). https://doi.org/10.1007/s10468-021-10083-5

\bibitem{gp} Gel'fand, I. M.; Ponomarev, V.A. \emph{Indecomposable representations of the Lorentz group}. (Russian). Uspehi Mat. Nauk 23(2(140)), 3–60. (1968).

\bibitem{hpr} Happel, D.; Preiser, U.; Ringel, C. M. \emph{Vinberg's chacterization of Dynkin diagrams using subadditive functions with application to $DTr$-periodic modules}. Proc. ICRA II (Ottawa, 1979), Lecture Notes in Math., vol. 832, Springer-Verlag, Berlin and New York, pp. 280-294. (1980).

\bibitem{ile} Ile, R. \emph{Change of rings in deformation theory of modules}. Trans. Amer. Math. Soc. 356, 4873-4896. (2004).

\bibitem{kra} Krause, H. \emph{Maps between tree and band modules}. J. Algebra \textbf{137.1}: 186-194. (1991).

\bibitem{kr21} Krause, H. \emph{Homological theory of representations}. (Vol. 195). Cambridge University Press. (2021).

\bibitem{lau} Laudal, O. A. \emph{Noncommutative deformations of modules.} in: The Roos Festschrift, vol. 2, Homology, Homotopy Appl. 4 357–396. (2002).

\bibitem{lop} L\'opez-Garc\'ia, D., Rizzo, P. and V\'elez-Marulanda, J.A. \emph{On weak universal deformation rings for objects of EXT-f\mbox inite categories of modules.} arXiv preprint, arXiv:2408.11259 (2024).

\bibitem{mgy} Magyar, A. \emph{Standard Koszul self-injective special biserial algebras.} J. Algebra Its Appl., Vol. 15, No. 03, 1650044. (2016). DOI: https://doi.org/10.1142/S0219498816500444.

\bibitem{msw} Meyer, D.; Soto, R.; Wackwitz, D. \emph{Universal deformation rings of modules for generalized Brauer tree algebras of polynomial growth.} Comm. Algebra 51, 3543 - 3555. (2023).

\bibitem{oz} Opper, S.; Zvonareva, A. \emph{Derived equivalence classif\mbox ication of Brauer graph algebras}. Adv. Math. 402, 108341. (2022). DOI:https://doi.org/10.1016/j.aim.2022.108341

\bibitem{jr} Rickard, J. \emph{Derived equivalences as derived functors}. J. Lond. Math. Soc. (2) 43, 37–48. (1991).

\bibitem{rv} Rizzo, P.; Vélez-Marulanda, J. \emph{A note on deformations of f\mbox inite dimensional modules over $\Bbbk$-algebras}. Comm. Algebra. {\bf 53.2}, 593-597. (2025). DOI: 10.1080/00927872.2024.2387048

\bibitem{sch} Schif\mbox f\mbox ler, R. \emph{Quiver representations}. Berlin: Springer. (2014).

\bibitem{schl} Schlessinger, M. \emph{Functors of Artin rings}. Trans. Amer. Math. Soc. 130, 208–222. (1968).

\bibitem{sib} Schroll, S. \emph{Trivial extensions of gentle algebras and Brauer graph algebras.} J. Algebra 444, 183–200. (2015). DOI: https://doi.org/10.1016/j.jalgebra.2015.07.037

\bibitem{syb} Schroll, S. \emph{Brauer graph algebras}. In: Homological methods, representation theory, and cluster algebras. pp. 177-223. (2018). DOI: 10.1007/978-3-319-74585-56

\bibitem{sk} Skartsaterhagen, O. \emph{Singular equivalence and the (Fg) condition.} J. Algebra 452, 66-93. (2016). DOI: https://doi.org/10.1016/j.jalgebra.2015.12.012

\bibitem{skw} Skowro\'nski, A. \emph{Periodicity in representation theory of algebras.} ICTP Lecture Notes, Advanced School and Conference on Representation Theory and Related Topics, 9–27 January 2006. (2006).

\bibitem{sy} Skowro\'nski, A.; Yamagata, K. \emph{Frobenius algebras. Vol. 1}. Eur. Math. Soc. (2011).

\bibitem{vm} Vélez-Marulanda, J. A. \emph{Universal deformation rings of strings modules over a certain symmetric special biserial algebra.} Beitr. Algebra Geom. 56(1):129-146. (2015). DOI: 10.1007/s13366-014-0201-y

\bibitem{yau} Yau, D. \emph{Deformation theory of modules}. Comm. Algebra 33, 2351–2359. (2005).

\end{thebibliography}
\end{document}